\documentclass[11pt,leqno]{article} 
\usepackage{amsmath,amssymb,latexsym,theorem,epsfig}
\usepackage[all]{xy}

\newcommand{\Appendix}[1]{%
  \refstepcounter{section}%
  \addcontentsline{toc}{section}%
    {\bfseries\appendixname~\thesection\ #1}%
    {\medskip\noindent \Large\bfseries\appendixname\ \thesection\ #1}%
\sectionmark{#1}\smallskip\noindent
\renewcommand{\theequation}{{\bf 
{{\thesection}.\arabic{subsection}}.{\arabic{equation}}}}
}

\makeatletter
\let\@internalcite\cite
\def\cite{\def\citeauthoryear##1##2{##1{##2}}\@internalcite}
\def\@biblabel#1{\def\citeauthoryear##1##2{##1{##2}}[#1]\hfill}
\makeatother

\theoremstyle{plain} 

\newtheorem{thm}{Theorem}
\newtheorem{theo}{Theorem}[section]
\newtheorem{lem}[theo]{Lemma}
\newtheorem{corr}[theo]{Corollary}
\newtheorem{con}[theo]{Conjecture}
\newtheorem{defi}[theo]{Definition}
\newtheorem{prop}[theo]{Proposition}
\newtheorem{prop-defi}[theo]{Proposition-Definition}
\newtheorem{lemma-defi}[theo]{Lemma-Definition}
{\theorembodyfont{\rmfamily} \newtheorem{rem}[theo]{Remark}}
{\theorembodyfont{\rmfamily} }
{\theorembodyfont{\rmfamily} \newtheorem{que}[theo]{Question}}

 \newcounter{inner}

\newcommand{\op}[1]{\operatorname{#1}}
\renewcommand{\theequation}{\thesection.\arabic{equation}}

\begin{document}

\title{Symplectic Lefschetz fibrations with arbitrary fundamental groups}
\author{J. Amor\'os \thanks{Partially supported by DGYCYT grant PB96-0234}, 
F. Bogomolov \thanks{Partially supported by NSF Grant DMS-9801591}, 
L. Katzarkov \thanks{Partially supported by NSF Grant DMS-
9700605 and A.P. Sloan research fellowship}, 
T. Pantev
\thanks{Partially supported by NSF Grant DMS-9800790} \\
(with an appendix by Ivan Smith)
}
\date{}
\maketitle

\begin{abstract}
In this paper we give an explicit construction of a symplectic
Lefschetz fibration whose total space is a smooth compact four
dimensional manifold with a prescribed fundamental group. We also
study the numerical properties of the sections in symplectic
Lefschetz fibrations and their relation to the structure of the
monodromy group.
\end{abstract}

\tableofcontents

\section{Introduction}
\setcounter{figure}{0}

In this paper we give an explicit construction of a symplectic
Lefschetz fibration whose total space is a smooth compact four
dimensional manifold with a prescribed fundamental group. 
The existence of such a fibration is also a consequence of the
remarkable recent 
work of Donaldson \cite{DON} (see also \cite{ARO}) who proved the 
existence of a Lefschetz
pencil  structure on any symplectic 4-manifold and the results of 
Gompf \cite{GOMPF}, who  proved that
any finitely presentable group can be realized as the fundamental 
group of a symplectic 4-manifold. 

Since Donaldson's proof is non-constructive, as an alternative 
we present a direct purely topological construction
of symplectic Lefschetz fibrations which is effective and allows an
explicit control on the number of singular fibers. The construction is
based on an algebraic geometric method for creating positive
relations among right handed Dehn twists. The ubiquity of such
relations combined with a simple group theoretic characterization 
of symplectic Lefschetz fibrations due to Gompf (see
Proposition~\ref{p:direct}) turns out to be sufficient for the
construction. 

Before we can state our main theorem we need to introduce some
notation. For any integer $n \geq 0$ denote by
$\pi_{n}$ the fundamental group of a compact Riemann surface of genus
$n$. 
As usual a group $G$ is called {\em finitely presentable} if it
can be written as a quotient of a free group on finitely many
generators by a subgroup generated by the conjugacy classes of
finitely many elements. By a {\em finite presentation} of a group $G$
we mean a surjective homomorphism $A \twoheadrightarrow G$ from some
finitely presentable group $A$ onto $G$ 
so that $\ker[A \to G]$ is generated as
a normal subgroup by finitely many elements in $A$. 

\begin{thm} Let $\Gamma$ be a finitely presentable group with a given 
finite presentation $a : \pi_g \twoheadrightarrow \Gamma$.
Then there exists a  surjective homomorphism  $b : \pi_h \to \pi_g$  
for some $h \geq g$ and
a symplectic Lefschetz fibration $f : X \to S^{2}$ such that
\begin{itemize}
\item[{\em (i)}] the regular fiber of $f$ is of genus $h$,
\item[{\em (ii)}] $\pi_{1}(X) \cong \Gamma$,
\item[{\em (iii)}] the natural
surjection of the fundamental group of the fiber of $f$ onto the 
fundamental group of $X$ coincides with $a\circ b$. 
\end{itemize}
\label{main}
\end{thm}

Note that any finitely presentable group $\Gamma$ admits a
finite presentation of the form $\pi_{g} \twoheadrightarrow \Gamma$ since
$\pi_{g}$ surjects onto a free group on $g$ generators.

\medskip

The map $b$ in the above theorem is not arbitrary. It factors as 
\[
\xymatrix{\pi_{h} \ar[rr]^{b} \ar[rd] & & \pi_{g} \\ & 
\pi_{e} \ar[ru]
}
\]
where the surface of genus $e$ is obtained form the surface of genus $g$
by adding handles and the surface of genus $h$ is obtained from
the surface of genus $e$ as a ramified finite covering.

\noindent
Our second theorem concerns symplectic fibrations of Lefschetz type
over curves of higher genus.

\begin{thm} \label{thm-high} 
Let $\Gamma$ be a finitely presentable group with a given 
presentation $a : \pi_g \twoheadrightarrow \Gamma$. Then there exist
a surjective homomorphism $\pi_e \to \pi_g$ 
and a symplectic Lefschetz fibration $f : X \to C_{k}$ over a smooth
surface $C_{k}$ of genus $k$
so that 
\begin{itemize}
\item[{\em (i)}] the regular fiber of $f$ is of genus $e$, 
\item[{\em (ii)}] $f$ has a unique singular fiber
\item[{\em (iii)}] the
fundamental group of $X$ fits in a short exact sequence
\[
1 \to \Gamma \to \pi_{1}(X) \to \pi_{k} \to 1,
\]
\item[{\em (iv)}] the natural surjective map from the fundamental group of
the fiber of $f$ to $\Gamma$  coincides with the composition $\pi_{e}
\to \pi_{g} \stackrel{a}{\to} \Gamma$.
\end{itemize}
\end{thm}

\

\medskip

This work is an elaboration on a discussion at the end of
\cite{BKS}. We describe in details an enhancement of the general  
technique for constructing  examples of symplectic fibrations used in 
\cite{BKS}. Our proof is based on exploiting the correspondence between 
subgroups of the mapping class group and graphs of vanishing cycles.

In general it is expected that every smooth four-dimensional manifold is
diffeomorphic to an achiral Lefschetz fibration possibly after some 
stabilization. The purpose of this paper is to study what other 
conditions besides chirality determine the SLF among all fibrations. 

The paper is organized as follows. In section two we give some 
preliminaries on subgroups of mapping class groups generated by Dehn
twists and recall an important result of Gompf which characterizes
symplectic Lefschetz fibrations via their monodromy representations. 
In section three  we explain the general construction and prove
Theorems \ref{main} and \ref{thm-high}. 
In section four  we give an explicit example of a symplectic Lefschetz
fibration of genus three  Riemann surfaces whose total space has 
first Betti number one
and a different  construction of a Lefschetz fibration whose total space
has a fundamental group isomorphic ${\mathbb Z}$. All this demonstrates
the flexibility of the construction for obtaining interesting
examples. 

The second
construction is similar in spirit to a computation done by Donaldson 
in which he has represented Thurston's example as a symplectic 
Lefschetz fibration of genus three Riemann surfaces.  A whole 
series of examples of the same flavor was
constructed independently in \cite{ozbagci-stipsicz} and
\cite{smithgenus2,smithoddgenus}. 

In the last section we apply the group theoretic part of the
construction to the study of the numerical properties of sections in
symplectic Lefschetz fibrations. Finally Appendix A, written by Ivan
Smith, presents a short proof of the non-existence of SLF with
monodromy contained in the Torelli group.

\bigskip
\bigskip

\noindent
{\bf Acknowledgements:} The authors would like to thank UC Irvine and
MPI-Bonn for hospitality. Parts of this work were done during the
authors participation in the UCI Hodge theory activities in June 1998
and during the first three
authors participation in the MPI activity dedicated to the work of 
B. Moishezon organized by F.Catanese, F.Hirzebruch and M.Teicher.
We would like  to thank the organizers of both events 
for the nice working atmosphere they have created. 

Special thanks are due to D. Auroux,
S.Donaldson, T.Fuller, R. Hain, B.Ozbagci, I.Smith, R.Stern and D.Toledo
for their nice suggestions and constant attention to this work. We also
thank D.Zagier and V.Platonov for useful group theoretic comments and
D. Auroux and R. Hain 
for pointing out some mistakes in a preliminary version of
this paper.

\subsection*{Notation and terminology}
\addcontentsline{toc}{subsection}{Notation and terminology}

\begin{description}
\item[$C_{g}$] a smooth compact oriented surface of genus $g$.
\item[$\Delta$] an analytic disk.
\item[$\Gamma$] a finitely presentable group.
\item[$\op{Map}_{g}(-)$] the subsemigroup in $\op{Map}_{g}$ generated
by the right Dehn twists $t_{r}, r \in R$.
\item[$\op{Map}_{g}^{n}$] the mapping  class group of a smooth genus 
$g$ surface with $n$ punctures.
\item[$\op{Map}_{g}{r}^{n}$] the mapping class group of a smooth genus
$g$ surface with $n$ punctures and $r$ boundary components.
\item[$\op{Map}_{R}$] the subgroup of the  mapping class group
$\op{Map}_{g}$  generated by the Dehn twists $t_{r}, r \in R$.
\item[$\op{Map}_{R}(-)$] the subsemigroup in $\op{Map}_{R}$ generated
by all conjugates of $t_{r}, r\in R$ within $\op{Map}_{R}$.
\item[$\op{Mon}$] the geometric monodromy group of a TLF, i.e. the  
image $\op{Mon} := \op{mon}(\pi_{1}(S^{2}\setminus\{q_{1},  \ldots,
q_{\mu} \},o))$ of the geometric monodromy representation.
\item[$\op{mon}$] the geometric monodromy representation associated
with a  TLF of genus $g$, i.e \linebreak $\op{mon} : \pi_{1}(S^{2} 
\setminus\{q_{1},\ldots, q_{\mu} \},o) \to \op{Map}_{g}$.
\item[$\mu$] 
the number of singular fibers in a TLF or a SLF.
\item[$o$] a base point in $S^{2}\setminus\{q_{1},  \ldots, q_{\mu} \}$.
\item[${\mathcal O}_{V}$] the structure sheaf  of an algebraic variety $V$.
\item[${\mathcal O}_{V}(1)$] a very  ample line bundle on an algebraic
variety $V$.
\item[$f : X \to S^{2}$] a topological or   symplectic Lefschetz
fibration (TLF or SLF).
\item[$f : X \to C_{k}$] a symplectic  fibration of Lefschetz type
over a  higher genus surface.
\item[$\pi_{g}^{n}$] the fundamental group of a  smooth genus $g$
surface  with $n$ punctures.
\item[$\pi_{n}$] the fundamental group of a smooth compete surface of
genus $n$.
\item[$\pi_{g} \twoheadrightarrow \Gamma$]
a finite presentation of $\Gamma$, i.e. a surjective homomorphism
whose  kernel is finitely generated as a normal subgroup.
\item[$Q_{i}$] a critical point of a topological or symplectic
Lefschetz  fibration.
\item[$q_{i}$] a critical value of a topological or symplectic
Lefschetz  fibration.
\item[$R \subset C_{g}$] a graph connected  collection of circles on
the  surface $C_{g}$.
\item[$s\subset C_{g}$] a circle in $C_{g}$ or in  other words a
smooth  connected one dimensional submanifold in $C_{g}$.
\item[$\Sigma$] a smooth projective algebraic curve.
\item[$T_{s}$] a right handed Dehn twist diffeomorphism associated
with a  circle $s \in C_{g}$. In other words for any $s \subset C_{g}$
one  chooses an orientation preserving identification  of a tubular  
neighborhood of $s$ with the oriented cylinder $[0,1]\times S^{1}
\subset  {\mathbb R}\times {\mathbb C}$ and then  defines $T_{s} \in  
\op{Diff}^{+}(C_{g})$ as the difeomorphism that acts a $(t,z) \mapsto 
(t,e^{2\pi it}z)$ on the cylinder and as identity everywhere else.
\item[$t_{s}$] the mapping class of a right handed Dehn twist $T_{s}$.
The  element $t_{s} \in \op{Map}_{g}$ depends only on the isotopy
class of  $T_{s}$.
\item[$U_{i}$] a small neighborhood of a critical value of a TLF or SLF.
\item[$U_{Q_{i}}$] a small neighborhood of a critical point in a TLF or SLF.
\item[$X_{i}$] the singular fiber of a TLF or a SLF  corresponding to
the  critical value $q_{i}$.
\item[$X_{o}$] a regular fiber of a TLF or a SLF.
\end{description}

\section{Symplectic Lefschetz fibrations}  
\setcounter{figure}{0}

First we recall some basic definitions and results. For more details
the reader may wish to consult \cite[Section~3.2.7]{sga7II} and \cite{kas}.

\begin{defi} Let $X$ be a smooth compact 4-manifold equipped with a 
smooth surjective map $f : X\to S^2$. We shall call it a {\em
topological Lefschetz fibration (TLF)} if the following conditions
hold: \label{defi-tlf}
\begin{list}{{\em (\roman{inner})}}{\usecounter{inner}}
\item The differential $df$ is surjective outside a finite subset of  
points $\{ Q_{1}, \ldots, Q_{\mu} \} \subset X$.
\item Whenever $p \in S^{2}\setminus \{ f(Q_{1}), \ldots, f(Q_{\mu})
\}$ the fiber $f^{-1}(p)$ is a  smooth orientable Riemann surface
of a given genus $g$.
\item The images $q_{i} := f(Q_i)$ are different for different $Q_i$.
\item Let $X_{i}$ denote the fiber of $f$ containing
$Q_{i}$. Then for any $i$ there are small disks $q_{i} \in
U_{i}$ and $Q_i \subset U_{Q_{i}} \subset X$ 
with $f : U_{Q_i} \to U_{i}$ being a complex Morse function in some complex
coordinates $(x,y)$ on $U_{Q_i}$ and $z$ on $U_{i}$, i.e. $z = f(x,y) =
x^{2} + y^{2}$. 
\end{list}
\end{defi}

\begin{figure}[!ht]
\begin{center}
\epsfig{file=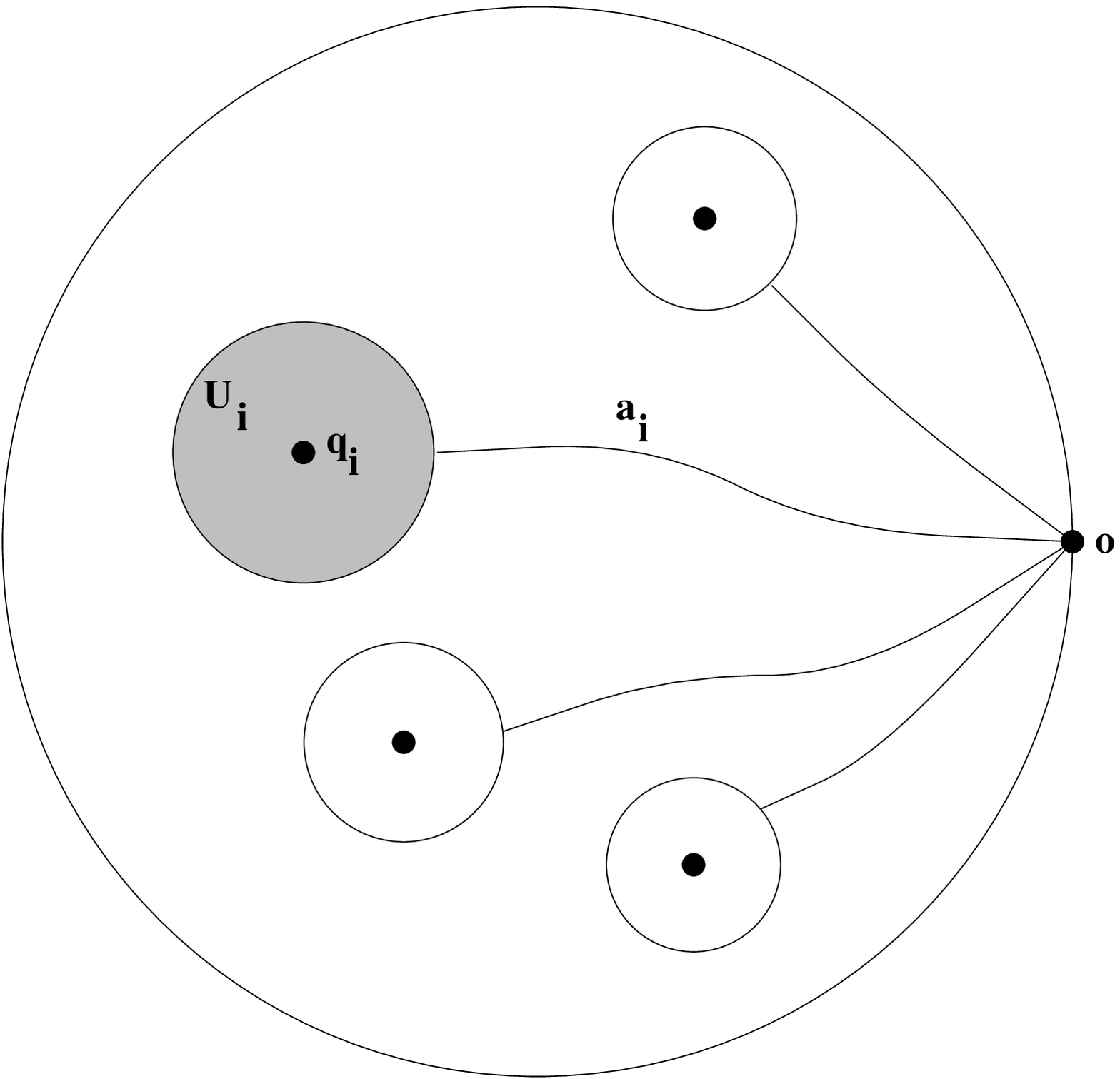,width=3in} 
\end{center}
\caption{An arc system for $f : X \to S^{2}$.}\label{fig1} 
\end{figure}

A topological Lefschetz fibration $f : X\to S^2$ will be called {\em
orientable} (or {\em chiral}) 
if there exists an orientation on $X$ so that the complex 
coordinates in (iv) above
can be chosen in a way compatible with the orientations on $X$ and
$S^{2}$. Note that the  definition of a topological Lefschetz
fibration is designed in such a way that the function $|f|^{2}$ is a
Morse function near the singular fibers $X_{i}$. The Morse flow
gives a handle body decomposition of $X$ and in particular one gets
standard retractions $cr_{i} : f^{-1}(\overline{U}_{i}) \to X_{i}$.
Fix a base point $o \in S^{2}$  which is a
regular value of $f$ and choose arcs $a_{1}, \ldots, a_{\mu}$  in
$S^{2}\setminus (\cup_{i} U_{i})$ which connect $o$ with some point on
$\partial U_{i}$, $i = 1, \ldots, \mu$  as in Figure~\ref{fig1}

\

\noindent
Such a collection of arcs and discs is called an {\em arc system} for
the TLF. A choice of an arc system gives a presentation of the
fundamental group of $S^{2}\setminus \{ q_{1}, \ldots, q_{\mu} \}$:
\[
\pi_{1}(S^{2}\setminus \{ q_{1}, \ldots, q_{\mu} \},o) = \langle
c_{1}, \ldots, c_{\mu} | c_{1}\cdot \ldots \cdot c_{\mu} = 1 \rangle,
\]
where geometrically $c_{i}$ is represented by the $o$-based loop in
$S^{2}\setminus \{ q_{1}, \ldots, q_{\mu} \}$ obtained by tracing
$a_{i}$ followed by tracing $\partial U_{i}$ counterclockwise and then
tracing back $a_{i}$ in the opposite direction.

Since the family $f : X \to S^{2}$ is locally trivial when restricted
to each $a_{i}$ we get well defined retractions of $X_{o}$ onto each
of the singular fibers $X_{i}$.  By abuse of notation these
retractions will be denoted by $cr_{i}$ as well. Each $cr_{i} : X_{o}
\to X_{i}$ contracts a smooth circle $s_i \subset X_{o}$  - {\em the
geometric vanishing cycle}. 
The boundary of $f^{-1}(U_{i})$ is diffeomorphic to a smooth
fiber bundle over the circle $\partial U_{i}$ with $X_{o}$ as a fiber. 
This fiber bundle is determined by a gluing diffeomorphism
$T_{s_{i}}: X_{o}\to X_{o}$ which is the usual Dehn twist along the circle
$s_i$.  More precisely $T_{s_{i}}$ is the right handed Dehn twist along
$s_{i}$ with respect to the orientation on $X_{o}$
compatible with the orientation on $U_{Q_{i}}$ given by the complex
coordinates around $Q_{i}$. The circles $s_{i}$ and the Dehn twists 
$T_{s_{i}}$ are uniquely determined up to a smooth
isotopy and thus give well defined elements $t_{i} \in
\op{Diff}^{+}(X_{o})/\op{Diff}^{+}_{0}(X_{o}) =: \op{Map}_{g}\subset
\op{Out}(\pi_{1}(X_{o}))$ - the group of mapping
classes of an oriented surface of genus $g$. The homomorphism 
$\op{mon}: \pi_{1}(S^{2}\setminus \{ q_{1}, \ldots, q_{\mu} \},o) \to
\op{Map}_{g}$, $\op{mon}(c_{i}) = t_{i}$ is called the {\em geometric
monodromy representation} of $f : X \to S^{2}$ and its image is called
the {\em geometric monodromy group} of the TLF. It is known
\cite[Theorem~2.4]{kas} that if 
$f : X \to S^{2}$ is an orientable TLF of genus $g \geq 2$, then the
geometric monodromy representation of $f$ uniquely determines 
the diffeomorphism type of $f$.

Note that if $C_{g}$ is an oriented surface and $s \subset C_{g}$ is a
smoothly embedded, homotopically non-trivial circle one can perform
both the right handed Dehn twist $T_{s} : C_{g} \to C_{g}$ and
the left handed Dehn twist $T_{s}^{-1} : C_{g} \to C_{g}$. If
however $f : X \to S^{2}$ is an orientable TLF and $C_{g} = X_{o}$
is given the induced orientation from $X$, then all of the geometric
monodromy transformations $\{T_{s_{1}}, \ldots, T_{s_{\mu}}
\}$ are right handed Dehn twists. This property actually characterizes
the orientable TLF completely.

We also consider topological Lefschetz fibrations which are compatible
with an additional closed non-degenerate $2$-form $w$ on $X$. 

\begin{defi} Assume that $f : X \to S^{2}$ is a TLF and that $(X,w)$
is a symplectic manifold. We say that $(f : X \to S^{2}; w)$ is a
{\em symplectic Lefschetz fibration (SLF)} if
for any $p\in S^2$ the form $w$ is non-degenerate on the fiber $X_p$
at $p$ in the sense that the smooth locus of $X_{p}$ is a symplectic 
submanifold in $X$ and for every $i$ the symplectic form $w_{Q_{i}}$ 
is non-degenerate on each of the two planes contained in the tangent cone of
$X_{i}$ at $Q_{i}$.
\end{defi}

Gompf \cite{gompf-stipsicz} 
had shown that under some mild restrictions the SLF 
can also be characterized in purely
topological terms. For the convenience of the reader we recall the
proof of this very useful fact. Different proofs can be found in
\cite{gompf-stipsicz} or \cite{smith}.

\begin{prop}[R.Gompf] \label{p:direct}
A topological Lefschetz fibration $f : X \to S^{2}$ of curves of genus
$g \geq 2$ admits a symplectic structure if and only if it is
orientable. 
\end{prop}
{\bf Proof.} Let us recall first the necessity of the orientation 
restriction. Assume that a symplectic Lefschetz fibration 
$f : X\rightarrow S^2$ has two singular values $q_1,q_2$ with singular
points $Q_1,Q_2$ above them, such that the monodromy
Dehn twists around $q_1,q_2$ have opposite orientations,
i.e. there exists a complex chart around $Q_1 \in X$ such
that $Q_1=(0,0)$, $f(x_{1},y_{1})= x_{1}^2+y_{1}^2$ and another chart
centered at $Q_2$ such that $f(x_{2},y_{2})=\bar{x}_{2}^2+ \bar{y}_{2}^2$. 
Here both charts are chosen to be compatible with the orientation on
$X$ given by the symplectic structure $w$. By the definition of a SLF
there is a symplectic form $w_{\op{base}}$ on $S^{2}$ induced from $w$
on $X$.

Join $q_1,q_2$ by a segment $I$ in the base $S^{2}$. The
fibration $f$ restricts to a trivial family over the interior
of $I$ with vanishing loop contractions at the 
endpoints.  Take a lift of the segment $I$ to $X$ 
with beginning and end points in the above 
coordinate charts, and choose a parallel trivial family of
horizontal tangent planes $\{ \pi_h(s) \}$ over this 
lifted segment. If we pick a trivial frame $e_1(s),e_2(s)$
for this family of planes we see that
$w_{\op{base}}(e_1(s),e_2(s))$ changes sign going from $q_1$
to $q_2$. This contradicts the continuity of the frame $e_1(s),e_2(s)$
and so the Dehn twists around $q_1,q_2$ must have the same
orientation.

We proceed to show that this condition is sufficient.

Let $f : X \to S^2$ be an orientable TLF with singular values $q_1,
\dots, q_{\mu}$. Take a topologically simple cover of $S^2$ by open
disks $\{ D_{\alpha} \}$, such that it includes a disk 
$U_i$ centered at every
singular value $q_i$, and $q_i \not\in \bar U_{\alpha}$ if $\alpha \neq
i$. We will put symplectic structures on the families over the disks
$U_{\alpha}$ first, and then glue them adapting an argument of Thurston
(see \cite[Theorem~6.3]{MS98}) for symplectic fibrations.

For every disk $D_i$ containing a singular value,
take the trivial family $C_g \times U_i \rightarrow U_i$,
endowed with a symplectic form $w_{\text{fiber}} \oplus
w_{\text{base}}$, the summands being symplectic forms on the factors.
By identifying $C_g$ with a regular fiber of $f_{| f^{-1}D_i}$ and
choosing an arc from the basepoint regular
value to the singular value of the fibration we get 
a vanishing loop $s_{i} \subset C_{g}$ that
determine the diffeomorphism type of the pencil $f$ over $D_i$.
We will perform symplectic surgery on the trivial family
$C_g \times U_i$ to make it diffeomorphic to $f$.

Let $q : {\mathbb C}^{2} \to {\mathbb C}$ be the standard quadratic
map $q(x,y) = x^{2} + y^{2}$.
Take a small ball $B$ centered at $(0,0) \in {\mathbb C}^2$,
a disk $D \subset q(B)$ such that the 3--sphere 
$\partial B$ is transverse to the fibers of $q$ and  such that $B =
q^{-1}(D)\cap \overline{B}$. The restricted family $q_{|B} : B \to D$
has a single quadratic singular fiber $B_{0}$ over $0 \in D \subset
{\mathbb C}$.

Let $s \subset C_g$ be one of the vanishing simple loops of 
the fibration $X$. Consider the normalization $\widetilde{B}_{0} \to
B_{0}$ and let $o_{1}$ and $o_{2}$ be the preimages of the singular point 
$(0,0) \in B_{0} \subset B$ in $\widetilde{B}_{0}$. Choose small
analytic discs in $\widetilde{B}_{0}$ centered at $o_{1}$ and $o_{2}$
respectively and let $\Delta_{1}$, $\Delta_{2}$ denote their images in
$B_{0}$.

\begin{figure}[!ht] 
\begin{center}
\epsfig{file=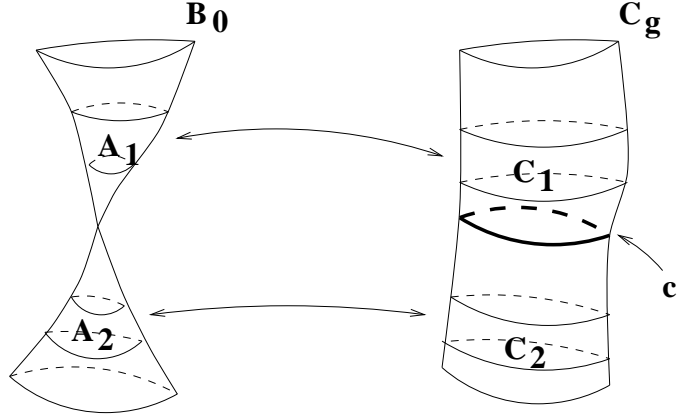,width=3.5in} 
\end{center}
\caption{The cylinders $A_{i}$ and $C_{i}$, $i = 1,2$.} \label{fig0}
\end{figure}

Select two open annuli $A_{1} \subset \Delta_{1}\setminus\{(0,0)\}$
and  $A_{2} \subset \Delta_{2}\setminus\{(0,0)\}$ such that 
the point $(0,0)$ does not lie in the closure $\bar A_1 \cup
\bar A_2$ (see Figure~\ref{fig0}). Using Moser's characterization 
of symplectic type of surfaces by volume, we may choose also two open
cylinders  $C_1,C_2$ on the opposite sides
of a bicollar neighborhood of $c$, such that both
$C_1,C_2$ are retracts of the bicollar neighborhood,
their adherence $\bar C_1 \cup \bar C_2$ does not intersect
$c$, and they are symplectomorphic to $A_1,A_2$ 
respectively. 

\bigskip

\noindent
The annuli $A_1,A_2,C_1,C_2$ are embedded in their
respective total spaces with trivial normal bundle,
so by Weinstein's symplectic
neighborhood theorem \cite[Lemma~2.1]{GOMPF},
theorem exists an $\varepsilon > 0$ and open neighborhoods
$W_i$ of $A_i$ in $B$, $V_j$ of $C_j$ in $C_g \times S^2$
for $i,j=1,2$ such that the $W_i, V_j$ are
symplectomorphic to $A_i \times D_{\varepsilon}, C_j \times
D_{\varepsilon}$ respectively. Shrinking $D$ if necessary,
we may now perform a
surgery by inserting the Dehn twist of $q_{|B} : B \rightarrow D$
by the identifications $W_i \cong V_i$. The symplectic
structures on $B$ and $C_g \times D_i$ define a symplectic structure on
$f$ over $D_i$ in this way.

We repeat this surgery over all the critical values of $f$,
rescaling the obtained symplectic structures $\omega_1, \dots,
\omega_{\mu}$ so that they induce the same symplectic structure $[
\sigma ]$ on
the regular fibers of the pencil. Through diffeomorphisms with trivial
families we may also endow the restrictions of $f$ over the
regular disks $U_{\alpha}$ with symplectic structures $\omega_{\alpha}$
inducing the same symplectic structure $[\sigma]$ on the regular fiber
of the pencil as the twist families over the singular values.

Let now $\tau_0 \in \Omega^2(X)$ be a closed 2--form such that its
restriction to the fibers represents the cohomology class $[\sigma]$.
Over every disk $D_{\alpha}$ we have that $\tau_0$ and the symplectic
form $\omega_{\alpha}$ are cohomologous, thus we may select 1--forms
$\lambda_{\alpha} \in \Omega^1(f^{-1}(D_{\alpha}))$ so that
$$
\omega_{\alpha}- \tau_0 = d \lambda_{\alpha}
$$
Choose now a partition of unity $\{ \rho_{\alpha} : D_{\alpha}
\rightarrow [0,1] \}$ subordinate to the cover $\{ U_{\alpha} \}$ and
such that for every critical value $q_i$ the function $\rho_i$ has
constant value $1$ on its neighbourhood. Define $\tau \in \Omega^2(X)$
by
$$
\tau= \tau_0+ \sum_{\alpha} d((\rho_{\alpha} \circ f) \lambda_{\alpha})
\, .
$$
This is a closed form, cohomologous to $\tau_0$, restricting to the
class $[\sigma]$ in every fiber and equal to the previously found
symplectic $\omega_i$ in neighbourhoods of the singular fibers. As in
the case of
smooth symplectic fibrations, this form $\tau$ is nondegenerate on the
tangent spaces to the fibers, and as it is defined on the compact total
space, if we select a symplectic form $\beta \in \Omega^2(S^2)$ the
forms
$$
\omega_K= \tau+ K f^* \beta
$$
are nondegenerate for sufficiently large $K$.

\  \hfill $\Box$

\section{The main construction}
\setcounter{figure}{0}

\subsection{Positive relations among right Dehn twists}
\label{subsec-positive} 

Let $g$ be  a nonnegative integer.
Fix a compact oriented reference surface $C_{g}$ of genus
$g$ and an infinite sequence of  distinct points $x_{0}, x_{1}, \ldots, x_{n},
\ldots  \in C_{g}$. 
Put $\pi_{g}^{n} :=  \pi_{1}(C_{g}\setminus \{x_{1}, \ldots, x_{n}
\},x_{0})$ and let $\op{Diff}^{+}(C_{g})^{n}_{r}$ denote the group of
all orientation preserving diffeomorphisms of $C_{g}$ that fix the
points $x_{1}, x_{2}, \ldots, x_{n+r}$ and induce the identity on the
tangent spaces $T_{x_{i}}C_{g}$ for $i = n+1, \ldots, n+r$. For any
triple of non-negative integers $(g,n,r)$ such that $2g - 2 + n + 2r >
0$ define the {\em mapping class group} $\op{Map}_{g,r}^{n}$ as the group
of connected components of $\op{Diff}^{+}(C_{g})^{n}_{r}$, i.e.
\[
\op{Map}_{g,r}^{n} := \pi_{0}(\op{Diff}^{+}(C_{g})^{n}_{r}).
\]
As usual we will skip the labels $n$ and $r$ if they happen to be
equal to zero. 

By definition the mapping class group $\op{Map}_{g}^{n}$ acts by
outer automorphisms on $\pi_{g}^{n}$. In fact this action identifies
\cite{looijenga-hain} 
$\op{Map}_{g}$ with the index two subgroup of $\op{Out}(\pi_{g})$
consisting of outer automorphisms acting trivially on
$H^{2}(\pi_{g},{\mathbb Z}) \cong H^{2}(C_{g},{\mathbb Z})$. Similarly
one can interpret the group $\op{Map}_{g}^{1}$ as the group of all
automorphisms of $\pi_{g}$ acting trivially on $H^{2}(\pi_{g},{\mathbb
Z})$. The group
$\op{Map}_{g,r}^{n}$ is generated by the right handed Dehn twists along all
(unoriented) non separating loops $c \subset C_{g}\setminus
\{x_{1}, \ldots, x_{n+r} \}$. 

\begin{rem} \label{remark-extension} 
For any $g,n$ the natural forgetful map $\op{Map}_{g,n} \to
\op{Map}_{g}^{n}$ is surjective and has a central kernel which can be
identified with the free abelian group generated by the Dehn twists
along simple loops around the punctures.
\end{rem}

For future reference define
$\op{Map}_{g,r}^{n}(-) \subset \op{Map}_{g,r}^{n}$ to be the
sub-semigroup of $\op{Map}_{g,r}^{n}$ generated by (the images of) 
all right handed Dehn twists. 

\medskip

The first step in the construction is the following simple observation.

\begin{lem} Let $s_{1}, \ldots, s_{m} \subset C_{g}$ be free simple closed
loops (not necessarily distinct)  
and let $t_{1}, \ldots, t_{m} \in \op{Map}_{g}$ be the
corresponding right handed Dehn twists. Suppose that there exist
integers $\{n_{i}\}_{i= 1}^{m}$ so that the $t_{i}$'s satisfy the
relation $\prod_{i}t_i^{n_i}=1$ in the mapping class group. 
Then there exists a TLF $f : X \to S^{2}$ with the following
properties:
\begin{list}{{\em (\alph{inner})}}{\usecounter{inner}} 
\item The regular fiber of $f$ can be identified with $C_{g}$ so that 
the vanishing cycles of $f$ are precisely the $s_{i}$'s.
\item The geometric monodromy group $\op{Mon}(f)$ of $f$ is just the
subgroup of the mapping class group generated by the $t_{i}$'s.
\item The fundamental group of $X$ is isomorphic to the quotient of
$\pi_{g}$ by the normal subgroup generated by all 
the $s_{i}$'s.
\item If all the $n_{i}$'s are positive, then $f$ is a SLF.
\end{list} \label{lem-tlf}
\end{lem}
{\bf Proof.} To construct the fibration $f : X \to S^{2}$ start with
the direct product $D_0\times C_{g}$  where $D_{0} \subset {\mathbb
C}$ is a small disk around zero. Next attach $|n_{i}|$-copies of 
small discs $D_i$, $i =
1, \ldots,  n$,  along the boundary of $D_{0}$. Over each $D_{i}$
choose a standard holomorphic Lefschetz fibration $U_{X_{i}} \to D_{i}$ 
with a unique singular fiber $X_{i}$ at
the center  so that the vanishing cycle is exactly $s_{i}$. 
The union of $D_0\times C_{g}$ and the fibrations $U(X_i)$ (each
appearing $|n_{i}|$-times respectively) has a structure of a topological
Lefschetz fibration $u : U \to D$ over a larger disc $D$. By
construction the Lefschetz fibration $u$ has a regular fiber
isomorphic to $C_{g}$ and vanishing cycles $s_{1}, \ldots,
s_{m}$. Moreover the 
geometric monodromy
transformation $\op{mon}(\partial D) \in \op{Map}_{g}$ for $u$ is precisely
the  product $\prod_{i}t_i^{n_i}$. Since
 the latter is equal to identity  in $\op{Map}_{g}$ we get that
$U_{|\partial D}$ is homotopy equivalent (and hence diffeomorphic) to a
product $C_{g}\times S^{1}$. In particular we can extend $u$ to a TLF
over $S^{2}$.

The conditions (a) and (b) are satisfied by construction. Let
$f^{\sharp} : X^{\sharp} \to S^{\sharp}$ denote the fibration obtained
from $f$ by removing the singular fibers. Since $f^{\sharp}$ is a
fiber bundle and by construction $f$ admits  a topological section we
can identify the fundamental group $\pi_{1}(X^{\sharp})$ with the
semidirect product $\pi_{1}(S^{\sharp})\ltimes_{\op{mon}} \pi_{g}$
where $\op{mon} : \pi_{1}(S^{\sharp}) \to \op{Aut}(\pi_{g})$ is the
monodromy representation. Now the Seifert-van Kampen theorem implies
that $\pi_{1}(X)$ is isomorphic to the quotient of $\pi_{g}$ by the
normal subgroup generated by the orbits of the vanishing loops $\{
s_{1}, \ldots, s_{m}\}$  under the monodromy group $\op{Mon}$ of $f$. On
the other hand any Dehn twist $t_{s}$ is uniquely characterized by the
property that the maximal quotient of $\pi_{g}$ on which $t_{s}$ acts
as the identity is the quotient of $\pi_{g}$ by the normal subgroup
generated by $s$. In particular for any $\gamma \in \pi_{g}$ the
element $\gamma^{-1}t_{s}(\gamma)$ is a product of conjugates of $s$
and so the normal subgroup of $\pi_{g}$ generated by the
$\op{Mon}$-orbits of $s_{1}, \ldots, s_{m}$ coincides with the normal
subgroup generated by $s_{1}, \ldots, s_{m}$ only which proves part
(c) of the Lemma.

Finally the fact that condition (d) holds for $f$ is a consequence of 
Proposition~\ref{p:direct}. The Lemma is proven.
\hfill $\Box$ 

\bigskip

The previous lemma shows that the  construction of symplectic Lefschetz
fibrations reduces to the problem of finding relations in the
semigroup $\op{Map}_{g}(-)  \subset \op{Map}_{g}$. In order to find
such relations and to be able to modify them we  will need to study
certain configurations of embedded circles in $C_{g}$.

\bigskip

\noindent
Recall \cite{ivanov-mccarthy} 
that given two isotopy classes $\rho$ and $\sigma$ of smooth
circles in $C_{g}$ one defines the {\em geometric intersection number}
$i(\rho,\sigma)$ as the minimum number of points of $r\cap s$ over all
representatives $r$ of $\rho$ and $s$ of $\sigma$. 
A finite set $R$ of smooth circles in $C_{g}$ is said to
be in {\em minimal position} if every
two elements $r, s \in R$ intersect transversally in exactly
$i([r],[s])$ points and no three elements in $R$ intersect.

Let $G_{R}$ be the dual graph of the one dimensional
cell complex $\cup_{r\in R} r \subset C_{g}$, i.e. the graph whose
vertices are the elements of $R$ and for which the edges connecting
two vertices $s, r \in R$ correspond to the intersection points of the
loops $s$ and $r$ in $C_{g}$.
 
\begin{defi} \begin{list}{{\em (\alph{inner})}}{\usecounter{inner}}
\item Two loops $s, r \in R$ will be called {\em
adjacent} if they intersect transversally at a single point, i.e. if
the vertices $s$ and $t$ in $G_{R}$ are connected by a single edge.
\item Two circles $s, r \in R$ will be called {\em graph connected} 
if the corresponding vertices are connected by a path of edges in
$G_{R}$. In other words $s, r \in R$ are graph connected if there 
exist a sequence of circles $s_{1}, \ldots, s_{m} \in
R$  so that $s = s_{1}$, $r = s_{m}$ and $s_{i}$ is adjacent to
$s_{i+1}$ for all $i = 1, \ldots, m-1$. 
\item Let $R$ and $S$ be two finite sets of circles in $C_{g}$ so that
$R\cup S$ is in minimal position. We will say that {\em $R$ is graph
connected to $S$} if  for any loop $r \in R$ there exists a 
loop $s \in S$ so that $r$ and $s$ are graph connected in $R \cup S$.
\end{list} \label{def-connected}
\end{defi}

\begin{rem} The adjacency condition in part (b) of
Definition~\ref{def-connected} is imposed only on two consecutive
loops in the sequence $s_{1}, \ldots, s_{m}$. In particular,
non-consecutive loops may intersect in an arbitrary transverse
fashion.
\end{rem}

\begin{figure}[!ht] 
\begin{center}
\epsfig{file=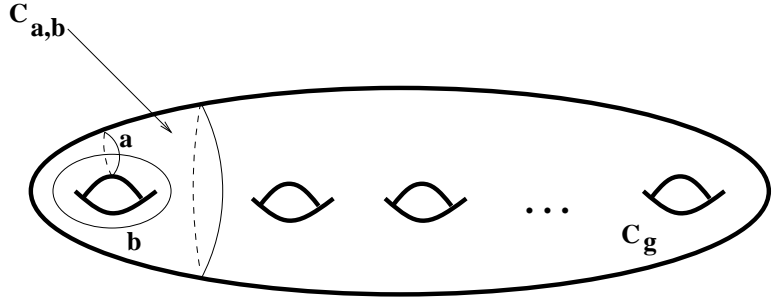,width=4in} 
\end{center}
\caption{A handle determined by two adjacent cycles.} \label{fig2}
\end{figure}

For any pair of smooth circles $a, b \subset C_{g}$ intersecting at
exactly one point one can find a neighborhood $C_{a,b}$ of $a\cup b$
in $C_{g}$ diffeomorphic to a torus with one hole as depicted on
Figure~\ref{fig2}. 
 
\bigskip

\noindent
Thus we obtain a
homomorphism $h_{a,b} : \op{Map}_{1,1} \to \op{Map}_{g}$ defined by
the handle $C_{a,b} \subset C_{g}$. 

For a set $R$ of circles in minimal position denote by $\op{Map}_{R}$
the subgroup of $\op{Map}_{g}$ generated by the right Dehn twists 
$\{ t_{r}\}_{r \in R}$. Let $\op{Map}_{R}(-)$ denote the
subsemigroup in $\op{Map}_{R}$ generated by the right twists $\{
t_{r}\}_{r \in R}$ and all of their conjugates in $\op{Map}_{R}$. Note
that $\op{Map}_{R}(-)$ is contained in $\op{Map}_{g}(-)$ as a subsemigroup 
since for every $\phi \in \op{Map}_{g}$ and every isotopy class $\rho$
of circles on $C_{g}$ one has $\phi\circ t_{\rho} \circ \phi^{-1} = 
t_{\phi(\rho)}$. Moreover by definition the semigroup $\op{Map}_{R}(-)
\subset \op{Map}_{R}$ is invariant under conjugation in $\op{Map}_{R}$.

The next two lemmas give a topological characterization  of the existence
of relations in $\op{Map}_{g}(-)$.

\begin{lem} Let $L$ be a finite set 
of circles in minimal position in $C_{g}$. Then $\op{Map}_{L}(-) =
\op{Map}_{L}$ if and only if there exists a finite relation $\alpha = 1$
where $\alpha$ is a product of only positive powers of $t_{\phi(l)}$'s
with $l \in L$, $\phi \in \op{Map}_{L}$ and
each $t_{l}, l \in L$ occurs at least once in $\alpha$. \label{lem-negative}
\end{lem}
{\bf Proof.} For the proof of the ``if'' part consider a relation $\alpha
= 1$ as in the statement of the lemma. By hypothesis there exist a 
positive integer $\mu$ and a
surjective map $l : \{1,2, \ldots, \mu \} \twoheadrightarrow L$ so
that 
\[
1 = \alpha = \prod_{i = 1}^{\mu} t_{l(i)}
\]
in $\op{Map}_{L}$. If we multiply both sides of the above relation by
$t_{l(1)}^{-1}$ on the left we get that 
\[
t_{l(1)}^{-1} = \prod_{i =2}^{\mu} t_{l(i)} \in \op{Map}_{L}(-).
\]
Since the relation is cyclic  this implies that $t_{l(i)}^{-1} \in
\op{Map}_{L}(-)$ for all $i = 1, \ldots, k$. Combining this with the
surjectivity of $l$ yields the inclusion $\{t_{l}\}_{l \in L} \subset
\op{Map}_{L}(-)$ and hence $\op{Map}_{L}(-) = \op{Map}_{L}$.

\medskip

To prove the ``only if'' part of the lemma note that the assumption
$\op{Map}_{L}(-) = \op{Map}_{L}$ implies that $t_{l}^{-1} \in
\op{Map}_{L}(-)$ for all $l \in L$. In other words each $t_{l}^{-1}$ can be
written as a product of finitely many of the right Dehn 
twists $\{ t_{\phi(s)} \}_{s \in L, \phi \in \op{Map}_{L}}$. 
Let now $m$ be the cardinality of $L$
and let $l_{1}, \ldots, l_{m}$ be an ordering of $L$. Next take 
the obvious relation
\[
1 = t_{l_{1}}t_{l_{1}}^{-1}t_{l_{2}}t_{l_{2}}^{-1}\ldots
t_{l_{m}}t_{l_{m}}^{-1}
\]
and then replace each of the left Dehn twists $t_{l_{i}}^{-1} \in
\op{Map}_{L} = \op{Map}(-)$ with the corresponding product of right
Dehn twists. The resulting right hand side will be a word $\alpha$
with the desired property. \hfill $\Box$

\bigskip

The next lemma gives a convenient criterion guaranteeing the existence
of positive relations among right Dehn twists.

\begin{lem} Let $R$ and $S$ be two sets of circles in minimal position
which are graph connected. 
Assume that $\op{Map}_{S}(-) =
\op{Map}_{S}$. Then $\op{Map}_{R \cup S}(-) = \op{Map}_{R \cup S}$.
\label{lem-connected}
\end{lem}
{\bf Proof.} By definition $\op{Map}_{R \cup S}(-)$ is the
subsemigroup in $\op{Map}_{g}$ generated by $t_{\phi(c)}$, $\phi \in
\op{Map}_{R \cup S}$, $c \in R\cup S$. Therefore we need to show that for any
$\phi \in \op{Map}_{R \cup S}$ and $c \in R$ the left twist
$t_{\phi(c)}^{-1}$ also belongs to $\op{Map}_{R \cup S}(-)$. Since
$t_{\phi(c)} = \phi\circ t_{c} \circ \phi^{-1}$ and $\op{Map}_{R \cup
S}(-)$ is conjugation invariant in $\op{Map}_{R \cup S}$ it suffices
to check that $t_{\phi(c)}^{-1} \in \op{Map}_{R \cup S}(-)$ whenever
$c \in R\cup S$.

If $c \in S$ then by assumption $t_{c}^{-1} \in \op{Map}_{S}(-)
\subset \op{Map}_{R\cup S}(-)$. If $c \in R$, then by the
hypothesis of the lemma $c$ is graph connected in $G_{R\cup S}$ with
some element $d \in S$. On the other hand according to
Lemma~\ref{lem-sl2} (b)  for any pair $a,b$ of adjacent
cycles on $C_{g}$ the right twists $t_{a}$, $t_{b}$ are
conjugate in $h_{a,b}(\op{Map}_{1,1}) \subset \op{Map}_{R\cup S}$. 
In particular 
$t_{c}^{-1}$ and $t_{d}^{-1}$ belong to the same conjugacy class in
$\op{Map}_{R\cup S}$. Since by assumption $t_{d}^{-1} \in \op{Map}_{S}(-)
\subset \op{Map}_{R\cup S}(-)$ and due to the conjugation invariance of
$\op{Map}_{R \cup S}(-)$ in $\op{Map}_{R \cup S}$ we conclude that 
$t_{c}^{-1} \in \op{Map}_{R \cup S}(-)$ as well, which concludes the
proof of the lemma. \hfill $\Box$

\bigskip

\noindent
Let now $\Gamma$ be a finitely presentable group and let $a : \pi_{g}
\twoheadrightarrow \Gamma$ be a given presentation.
As we can see from the previous lemmas the construction of a SLF whose
total space has fundamental group $\Gamma$ reduces to finding a
graph connected system of circles $R$ so that $R$ generates $\ker(a)$
as a normal subgroup and the right Dehn twists about
some subsystem of $S$ satisfy a positive relation in
$\op{Map}_{g}$. 

To achieve this we will have to modify the presentation $a$ and in
particular enlarge the genus $g$. This will be done in two steps which
are explained in the next two sections. 

\subsection{Geometric presentations}

The first step is purely topological. Starting with a presentation $a : \pi_{g}
\twoheadrightarrow \Gamma$ we show how to add handles to $C_{g}$ to
obtain a new presentation $\psi : \pi_{e} \twoheadrightarrow \Gamma$
for which the cycles generating $\ker(\psi)$ are nicely situated on
$C_{e}$.

We begin with the following definition:

\begin{defi} Let $\Gamma$ be a finitely presentable group. A {\em
geometric presentation} of $\Gamma$ is a  
surjective homomorphism $\psi : \pi_{e} \to \Gamma$, where $\pi_e$ is 
the fundamental group of a compact Riemann surface $C_{e}$ of genus
$e$ such that: 
\begin{list}{{\em (\alph{inner})}}{\usecounter{inner}}
\item The group of relations $\ker(a)$ is generated as a
normal subgroup by finitely many simple closed loops $r_{1},
\ldots, r_{m} \subset C_{e}$ ;
\item Every two of the loops $r_{1}, \ldots, r_{m}$ intersect
transversally at most at one point and the subspace $\cup_{i}r_{i}
\subset C_{e}$ is connected.
\end{list} \label{defi-geometric-presentation}
\end{defi}

As we will see in the next section the geometrically presented groups
are well suited for algebraic geometric manipulations. Thus it is
important to find a procedure for constructing geometric presentations
of finitely presentable groups. We have the following 

\begin{lem} \label{lem-geometric-presentation} Let $a : \pi_{g}
\twoheadrightarrow \Gamma$ be a given presentation. Then there exists
a surjective homomorphism $\pi_{e} \twoheadrightarrow \pi_{g}$ so that the
composition $\pi_{e} \to \pi_{g} \to \Gamma$ is a geometric presentation.
\end{lem}
{\bf Proof.} Let $R \subset \ker(a)$ be a finite subset of elements
which generates $\ker(a)$ as a normal subgroup (such a set exists
since by assumption $a$ is a finite presentation of
$\Gamma$). Represent each element $r \in R$ by a free immersed loop
$c_{r} \subset C_{g}$ so that the one dimensional simplicial
sub complex $M := \cup_{r \in R} c_{r} \subset C_{g}$ is connected and has
only ordinary double points as singularities. Such a collection of
immersed loops $c_{r}$ can be easily found as a small perturbation of a
standard representation of the elements in $R$ as immersed loops based
at some fixed point $x_{0} \in C_{g}$.

\begin{figure}[!ht] 
\begin{center}
\epsfig{file=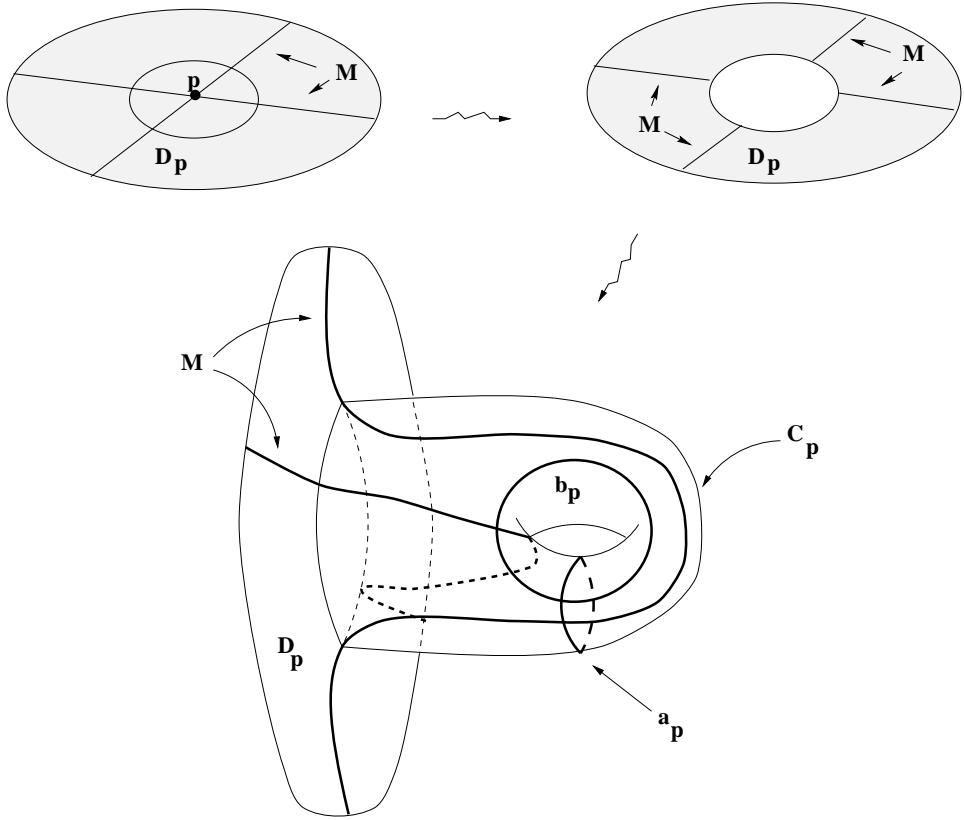,width=5in} 
\end{center}
\caption{Gluing the handle $C_{p}$ to $D_{p}\setminus \{p\}$.} \label{fig3}
\end{figure}

For each singular point $p \in M$ choose a
small disk $p \in D_{p} \subset C_{g}$ which doesn't contain any other
singularity of $M$. Now for each $p \in \op{Sing}(M)$ delete from
$D_{p}$ a smaller disk centered at $p$ and glue a handle $C_{p}$ to
$C_{g}$ along the inner rim of the resulting annulus. The two branches
of $M\cap D_{p}$ meeting at $p$ can be then completed to two smooth
disjoint curves on $C_{p}$ as shown on Figure~\ref{fig3}.

\

\noindent
In this way we obtain a new smooth surface 
\[
C_{e} := (C_{g}\setminus \op{Sing}(M))\cup
(\cup_{p \in \op{Sing}(M)} C_{p})
\]
of genus $e = g + \# (\op{Sing}(M))$ and a system of smooth disjoint circles 
$\tilde{c}_{r} \subset C_{e}$, $r \in R$. 
Next for every $p \in
\op{Sing}(M)$ we choose standard generators $a_{p}, b_{p} \subset C_{p}
\subset C_{e}$ of the fundamental group of $C_{p}$ as on
Figure~\ref{fig3}. 
By construction $\Gamma$ is isomorphic to the quotient  
of $\pi_{e}$ by the normal subgroup $N \triangleleft \pi_{e}$
generated by the finite set of elements 
\[
\{ \tilde{c}_{r}\}_{r \in
R}\cup (\cup_{p \in \op{Sing}(M)}\{ a_{p}, b_{p} \} \subset \pi_{e}.
\]
But we have glued the handles in such a way that the $\tilde{c}_{r}$'s
are disjoint from each other, $(a_{p}\cup b_{p})\cap (a_{q}\cup b_{q})
\neq \varnothing$ only if $p=q$ and each $a_{p}$ or $b_{p}$ intersects
exactly one of the $\tilde{c}_{r}$'s transversally at a single
point. Finally since every $a_{p}$ intersects $b_{p}$ at one point we
conclude that the union of all these cycles is connected in $C_{e}$
and hence $\pi_{e} \to \pi_{e}/N \cong \Gamma$ is a
geometric presentation. \hfill $\Box$

\subsection{The proof of Theorem~\ref{main}}
Assume now that $\Gamma$ is a group with a fixed
geometric presentation $\psi: \pi_{e} \to \Gamma$. Let $R = \{r_{1},
\ldots, r_{m} \}$ be a set of circles in minimal position in $C_{e}$
so that $R$ generates $\ker(a)$ as a normal subgroup as in
Definition~\ref{defi-geometric-presentation}. Note that
without a loss of generality we may assume that $R$ contains a 
non-separating circle $s \in R$. Indeed, if $\ker(a)$ happens to be
contained in $[\pi_{e},\pi_{e}]$ we can always glue an extra handle to
$C_{e}$ and add to $R$ the three standard generators of the first
homology of the handle (plus some extra loops if required).

Let $\Sigma$ be a smooth projective algebraic curve of genus $e$ which
we have identified as a $C^{\infty}$ manifold with $C_{e}$.  Fix a
point $p \in s \subset \Sigma$ which does not lie on any of the
circles in $R\setminus \{s\}$ 
and let $p \in \Delta \subset \Sigma$ be a small analytic disc which
is disjoint from all the circles in $R\setminus \{s\}$.

Let $V = \Sigma\times {\mathbb P}^1$ and let $D \subset
V$ be a very ample divisor such that $D = \sum_{i} D_{i}$ with each
$D_{i} \subset V$ being a section for the projection $p_{\Sigma} : V
\to \Sigma$. By replacing ${\mathcal O}_{V}(D)$ by its third power if
necessary we may further assume that the degree of $D$ on 
${\mathbb P}^{1}$ is divisible by three.

Let $p_{i}$, $\Delta_{i}$ and $s_{i}$ denote the 
preimages  in $D_{i}$ of $p$, $\Delta$ and $s$ respectively. Similarly
let $R_{i}$ be the set of circles in $D_{i}$ consisting of the 
preimages of the circles in $R$ via the projection $p_{\Sigma|D_{i}} :
D_{i} \to \Sigma$. When
${\mathcal O}_{V}(D)$ is chosen to be sufficiently ample and the
divisor $D$ is chosen to be a generic deformation of a set of sections
of $p_{\Sigma}$ which pass trough a given point lying over $p \in
\Sigma$ we may easily arrange that all intersections of the $D_{i}$'s 
are transverse and that for every pair of indices $i\neq j$ we have 
$\Delta_{i}\cap \Delta_{j} \neq \varnothing$ (see Figure~\ref{fig4} below). 

\begin{figure}[!ht] 
\begin{center}
\epsfig{file=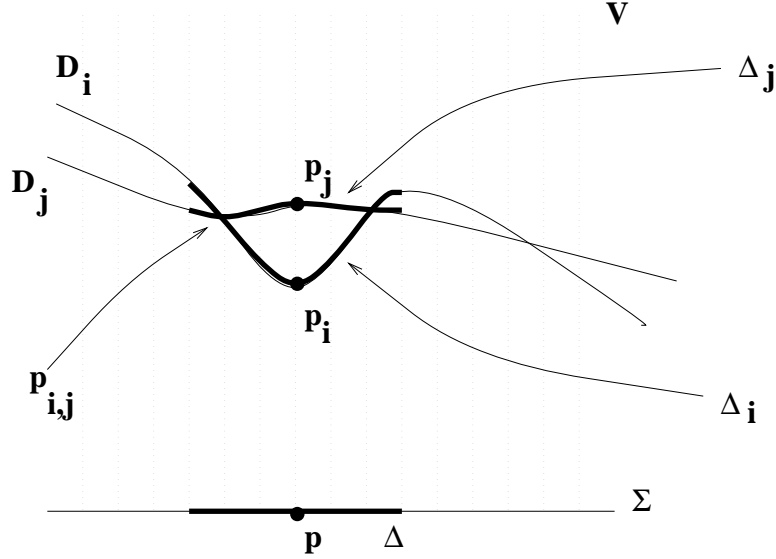,width=4in} 
\end{center}
\caption{The divisors $D_{i}$.} \label{fig4}
\end{figure}

\

\begin{figure}[!ht] 
\begin{center}
\epsfig{file=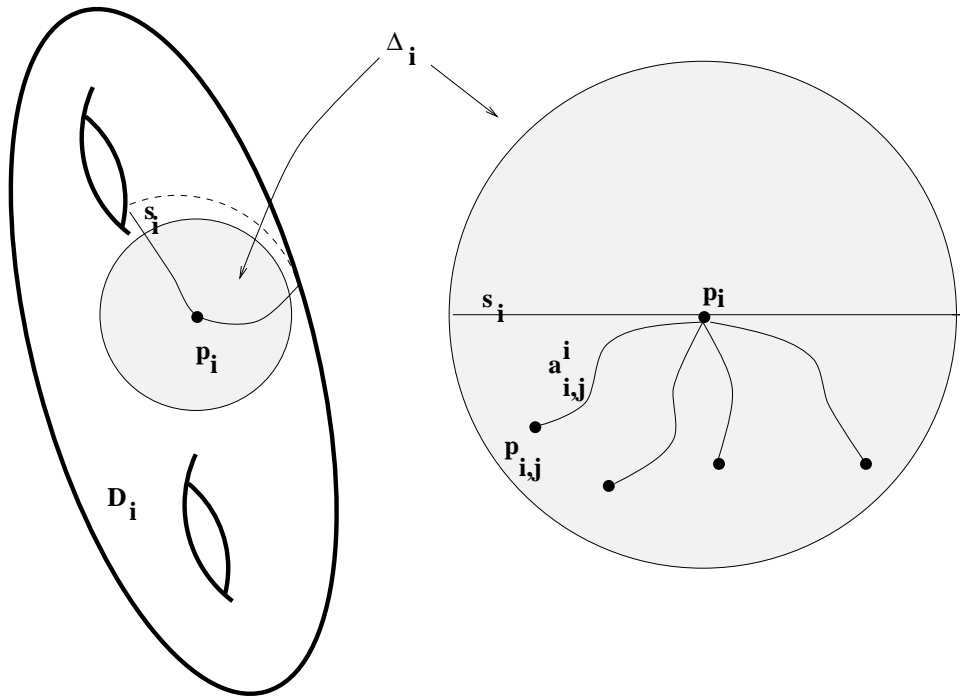,width=5in} 
\end{center}
\caption{The system of arcs in $\Delta_{i}$.} \label{fig5}
\end{figure}

\noindent
For any $i \neq j$ pick a point $p_{i,j} \in  \Delta_{i}\cap \Delta_{j}$. 
Choose arcs $a^i_{i,j} \subset \Delta_{i}$ connecting $p_{i}$ with
$p_{i,j}$ which meet only at $p_{i}$ and do not intersect $s_{i}$ at
any other point (see Figure~\ref{fig5}).

\begin{rem} For each $i$ and $j$ only {\em one} point is chosen in
$\Delta_{i}\cap \Delta_{j}$. Therefore the set of points $\{p_{i,j}\}$
will be only a subset in $\cup_{i<j}(D_{i}\cap D_{j})$ in general.
\end{rem}

\

\bigskip

\noindent
Consider a generic pencil in the linear system $|D|$ which contains $D$
as a member and has a smooth general member. After blowing up the base
points of this pencil on $V$ we obtain a smooth surface $\widehat{V}$
and a projective morphism $v : \widehat{V} \to {\mathbb C}{\mathbb
P}^{1}$ having $D$ as a fiber over some point $d \in {\mathbb C}{\mathbb
P}^{1}$ and a smooth connected
general fiber. Let $V_{o}$ be a smooth fiber of $v$ over some
reference point $o \in {\mathbb C}{\mathbb P}^{1}$ which is close to
$d$. Then we have (once we choose an arc system for $v$) a well
defined deformation retraction $cr : V_{o} \to D$ which collapses
certain smooth circles on $V_{o}$ to the points in the finite set
$\cup_{i<j}(D_{i}\cap D_{j})$. Note that trough the retraction $cr$ we
may view the $\Delta_{i}\setminus \{p_{i,j}\}_{j\neq i}$ as 
punctured discs on the smooth curve $V_{o}$, the set
$\cup_{i} R_{i}$ as a set of circles in minimal position on
$V_{o}$ and the arcs $r^{i}_{i,j}$ as arcs on $V_{o}$. Moreover
by a slight perturbation within each $\Delta_{i}$ we we can arrange
that for each pair of indices $i \neq j$ the arcs $r^{i}_{i,j},
r^{j}_{i,j}  \in V_{o}$ in joins smoothly at a point on the circle
$c_{i,j}$. 

\begin{figure}[!ht] 
\begin{center}
\epsfig{file=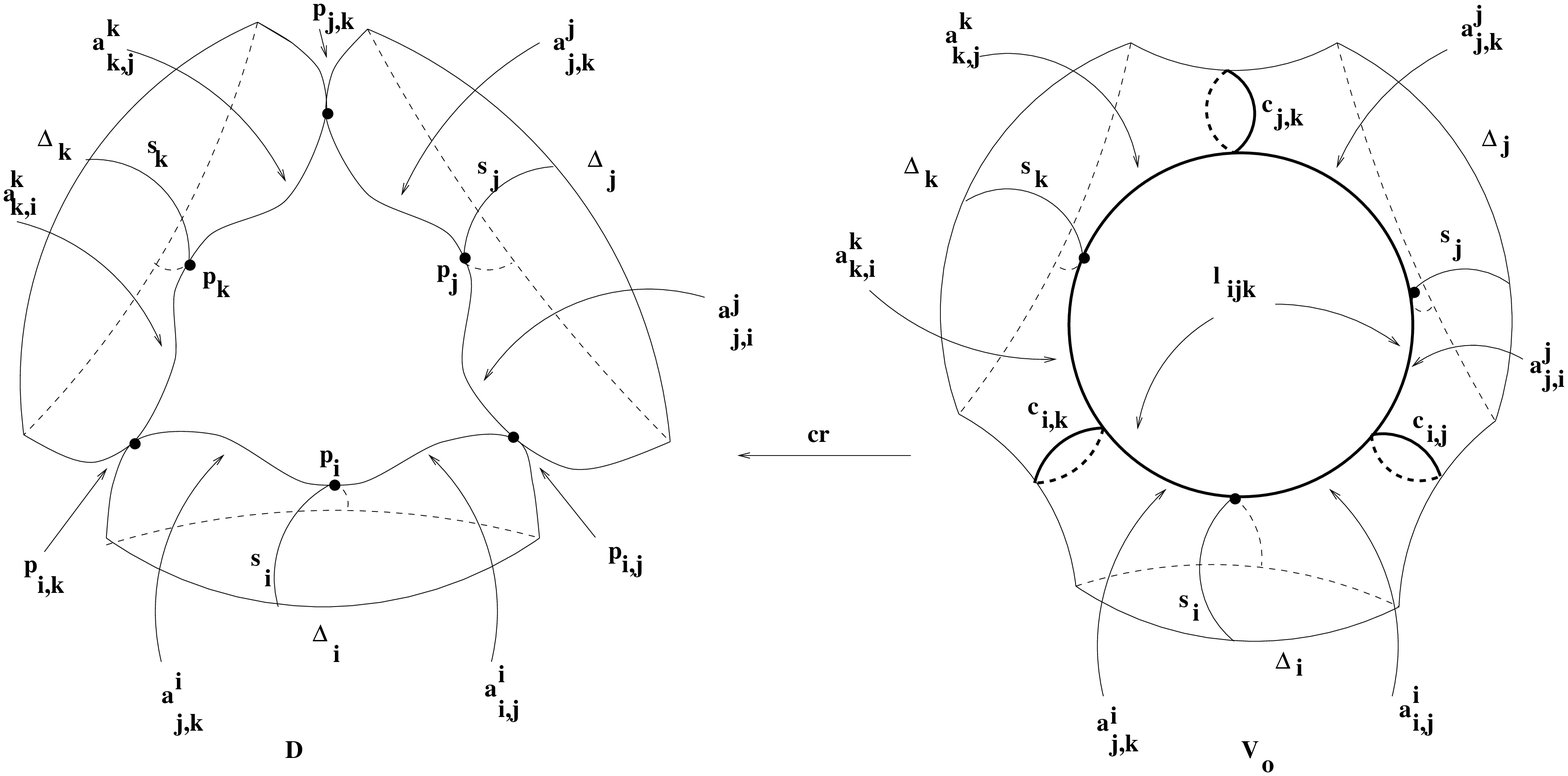,width=6.4in}
\end{center}
\caption{Deformation retraction of $V_{o}$ onto $D$.}  \label{fig6}
\end{figure}

In particular we see that the points $p_{i,j} \in D$ correspond
to vanishing cycles $c_{i,j} \subset V_{o}$ for the pencil $v$ which 
are disjoint from $s_i \subset V_{o}$. Moreover
by a slight perturbation within each $\Delta_{i}$ we can arrange
that for each pair of indices $i \neq j$ the arcs $r^{i}_{i,j},
r^{j}_{i,j}  \in V_{o}$ in joins smoothly at a point on the circle
$c_{i,j}$ (see Figure~\ref{fig6}).

\

\noindent
For each triple of distinct indices $i,j,k$ denote by $l_{i,j,k}$ the
smooth unoriented circle in $V_{o}$ obtained by tracing 
the segments $a^{x}_{x,y}$,
$\{x, y\} \subset \{i,j,k\}$ in the following order (see Figure~\ref{fig6})
\[
l_{i,j,k} = a^{k}_{k,i}\circ \bar{a}^{k}_{k,j}\circ a^{j}_{j,k}\circ
\bar{a}^{j}_{j,i}\circ a^{i}_{i,j}\circ \bar{a}^{i}_{i,k},
\]
where we have adopted the convention, that $a^{x}_{xy}$ is oriented
from $p_{x}$ to $p_{x,y}$ and $\bar{a}^{x}_{xy}$ is oriented from
$p_{x,y}$ to $p_{x}$.

By construction, the $l_{ijk}$'s are free loops representing
elements in the kernel of the natural surjection $(p_{\Sigma|V_{o}})_{*} :
\pi_{1}(V_{o}) \twoheadrightarrow \pi_{1}(\Sigma)$.  
 
Before we state the main result of this section we need to introduce
some notation. Let $h$ be the genus of $V_{o}$ and let $S
\subset \pi_{1}(V_{o})$ be the set of geometric
vanishing cycles for the algebraic pencil $v$. We have the following:

\begin{prop} Let $L$ be the set of circles in minimal position on $C_{h}$
defined as
\[
L := (\cup_{i} R_{i})\cup S \cup \{ l_{123}, l_{456}, \ldots
\}.  
\]
Then there exists a SLF $f : X \to S^{2}$ of genus $h$ such that 
\begin{itemize}
\item The set of vanishing cycles of $f$ is exactly $L$.
\item The fundamental group of $X$ is canonically isomorphic to
$\Gamma$. 
\item The identifications $\pi_{1}(X) \cong \Gamma$ and $V_{o} \cong
C_{h}$ can be chosen so that the natural epimorphism  
$\pi_{h} \to \pi_{1}(X)$ (induced from the inclusion of $C_{h}$ in $X$
as a regular fiber of $f$) becomes the composition 
\[
\xymatrix@1{\pi_{h} = \pi_{1}(V_{o}) \ar[rr]^-{(p_{\Sigma|V_{o}})_{*}}
& &
\pi_{1}(\Sigma) = \pi_{e} \ar[r]^-{\psi} & \Gamma
}
\]
\end{itemize} \label{prop-main}
\end{prop}
{\bf Proof.} Let $\iota : V_{o} \hookrightarrow V$ be the natural
inclusion of the divisor $V_{o}$ in the surface $V$. Consider the
induced map $\iota_{*} : \pi_{1}(V_{o}) \to \pi_{1}(X)$. Since $V_{o}$
is the regular fiber of the algebraic Lefschetz fibration $v :
\widehat{V} \to {\mathbb P}^{1}$ we conclude by Lemma~\ref{lem-tlf} (c)
that $\ker(\iota_{*})$ is the normal subgroup generated by the circles
in $S$. On the other hand since $p_{\Sigma *}$ identifies
$\pi_{1}(V)$ with $\pi_{1}(\Sigma)$ we have that $\ker(\iota_{*}) =
\ker((p_{\Sigma|V_{o}})_{*})$ and so $l_{ijk} \in \ker(\iota_{*})$ for
all triples of indices $i, j, k$. Due to this observation the normal
subgroup in $\pi_{1}(V_{o}) = \pi_{h}$ generated by the elements in
$L$ actually coincides with the normal subgroup generated by the
elements in $(\cup_{i} R_{i}) \cup S$. Moreover by the handle body
decomposition for the Lefschetz pencil $v$ we have a commutative
diagram
\[
\xymatrix{\pi_{1}(V_{o}) \ar[rr]^-{(p_{\Sigma|V_{o}})_{*}}
\ar[dr]_{cr_{*}} & & \pi_{1}(\Sigma) \\
& \pi_{1}(D) \ar[ur]_-{(p_{\Sigma|D})_{*}} &
}
\]
and hence each set $R_{i} \subset \pi_{1}(V_{o}) = \pi_{h}$ is mapped 
bijectively to the set $R \subset \pi_{e}$ by the map $\iota_{*}$. 
This shows that $L \subset \ker(\psi\circ\iota_{*})$ and that moreover
$L$ generates $\ker(\psi\circ\iota_{*})$ as a normal subgroup in
$\pi_{h}$. Therefore By Lemma~\ref{lem-tlf} the proposition will be
proven if we can find a positive relation among the elements in $L$.

The set of circles $S$ is the set of vanishing cycles in an algebraic
Lefschetz pencil and so we have a positive relation among the elements
in $S$ as in Lemma~\ref{lem-tlf}. Thus $\op{Map}_{S}(-) = \op{Map}_{S}$ by 
the ``if'' part of Lemma~\ref{lem-negative}.

Put $R := \cup_{i} R_{i} \cup \{ l_{123}, l_{456}, \ldots \}$. By
the geometric presentation assumption we know that each $R_{i}$ 
is graph connected. Also by construction
$l_{ijk}$ is graph connected to the circles $s_{i} \in R_{i}$, $s_{j}
\in R_{j}$ and $s_{k} \in R_{k}$ and to the vanishing cycles $c_{i,j},
c_{j,k}, c_{i,k} \in S$. Therefore $R$ is graph connected to
$S$ and so by Lemma~\ref{lem-connected} we have $\op{Map}_{R\cup S}(-)
= \op{Map}_{R\cup S}$. Now we can apply the ``only if'' part of 
Lemma~\ref{lem-negative} to the set $L = R\cup S$ which concludes the
proof of the proposition. \hfill $\Box$

\bigskip

Granted the previous proposition Theorem~\ref{main} now becomes a
tautology: 

\bigskip

\noindent
{\bf Proof of Theorem~\ref{main}.} Given a finite presentation $a :
\pi_{g} \to \Gamma$ of $\Gamma$ first use
Lemma~\ref{lem-geometric-presentation} to obtain a geometric
presentation $\psi : \pi_{e} \to \Gamma$ and then apply
Proposition~\ref{prop-main} to obtain a SLF with fundamental group
$\Gamma$. \hfill $\Box$

\subsection{Symplectic fibrations over bases of higher genus}

The construction explained in the two previous sections can be easily
modified to produce examples of symplectic Lefschetz fibrations over
surfaces of genus bigger than zero (see Remark~\ref{rem-high} for an
explanation of the terminology). Moreover in that case
the algebraic geometric part of the construction becomes superfluous 
and the resulting fibrations have monodromy groups generated
essentially by the Dehn twists with respect to the circles giving the
relations in some geometric presentation. This phenomenon is
illustrated in Theorem~\ref{thm-high} which we prove below. 

\begin{rem} \label{rem-high}
Before we prove the theorem a terminological remark is in
order. Note that the notion of a symplectic Lefschetz fibration introduced in
Definition~\ref{defi-tlf} makes sense without the assumptions that the
base of the fibration is a sphere and that the total space is four
dimensional. In fact the whole setup of Definition~\ref{defi-tlf}
easily generalizes to fibrations of even dimensional manifolds over a
base of even dimension.

In view of this we will adopt a slightly different terminology for the
remainder of this section. Namely a symplectic Lefschetz fibration
will mean a map $f : X \to C$ from a smooth 4-fold $X$ to a smooth
compact Riemann surface $C$ which satisfies all the conditions (i)-(iv)
from Definition~\ref{defi-tlf} but with $S^{2}$ replaced by $C$. 

This abuse of terminology should not cause any confusion since
symplectic Lefschetz fibrations like that will be considered only in
this section.
\end{rem}

We can now state the main result of this section

\bigskip

\begin{prop} \label{prop-high}
Let $\Gamma$ be a finitely presentable group with a given 
presentation $a : \pi_g \twoheadrightarrow \Gamma$. Then there exist
a geometric presentation $\pi_e \to \pi_g \stackrel{a}{\to} \Gamma$ 
and a symplectic Lefschetz fibration $f : X \to C_{k}$ for some $k >>
0$ so that 
\begin{itemize}
\item the regular fiber of $f$ is of genus $e$, 
\item the
fundamental group of $X$ fits in a short exact sequence
\[
1 \to \Gamma \to \pi_{1}(X) \to \pi_{k} \to 1
\]
\item the natural surjective map from the fundamental group of
the fiber of $f$ to $\Gamma$  coincides with the composition $\pi_{e}
\to \pi_{g} \stackrel{a}{\to} \Gamma$.
\end{itemize}
\end{prop}
{\bf Proof.} 
Let $\psi : \pi_{d} \to \pi_{g} \to \Gamma$ be a geometric
presentation as in Lemma~\ref{lem-geometric-presentation}. Let $R$ be
a generating set of circles for $\ker(\psi)$ as in
Definition~\ref{defi-geometric-presentation}. Even though the circles
in $R$ need not satisfy a positive relation in $\op{Map}_{d}$ the
subsemigroup $\op{Map}_{R}(-)$ is very close to being the whole group
$\op{Map}_{R}$ as the following arguments show.

\begin{lem} The abelianization of $\op{Map}_{R}$ is a cyclic
group. \label{lem-cyclic}
\end{lem}
{\bf Proof.} Fix a reference element $c \in R$. 
The group $\op{Map}_{R}$ is generated by the Dehn twists 
$\{t_{s} \}_{s \in R}$. Since by definition any two circles in $R$ are
graph connected we conclude as in the proof of
Lemma~\ref{lem-connected} that all generators of $\op{Map}_{R}$ are
conjugate to each other.  Therefore any character  $\op{Map}_{R} \to
S^{1}$ is defined by its value on $t_c$ and so the quotient 
$\op{Map}_{R}/[\op{Map}_{R},\op{Map}_{R}]$ is generated by
$t_{c}[\op{Map}_{R},\op{Map}_{R}]$. The lemma is proven. \hfill $\Box$

\medskip

\noindent
The previous lemma shows that if we can find an integer $n$ 
such that $t^n_s  \in
[\op{Map}_{R},\op{Map}_{R}]$ for any $s\in R$, then the abelianization
of $\op{Map}_{R}$ is of order at most $n$. In fact the same
statement holds for any subgroup $G$ of the mapping class group
$G\subset \op{Map}_{d}$ which is contained in 
the normal closure of $\op{Map}_{R}$ in $\op{Map}_{d}$.

\medskip

\begin{figure}[!ht]
\begin{center}
\epsfig{file=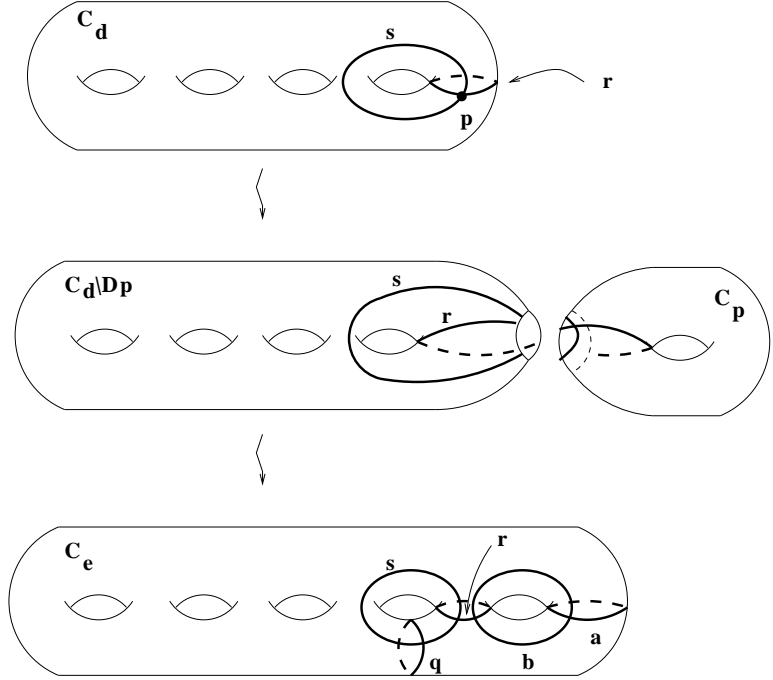,width=4in}
\end{center}
\caption{Enlarging a geometric presentation.}  \label{fig7}
\end{figure}

For the next step in the proof we will need to assume that the set $R$
is big enough. More precisely let $r, s \in R$ be two adjacent circles
on $C_{d}$ and let $C_{r,s} \subset C_{d}$ be the corresponding
handle. Let $p$ be the intersection point of $r$ and $s$.
By deleting a small disk $D_{p} \subset C_{d}$ centered at $p$
and gluing a  handle $C_{p}$ in its place as in the proof of
Lemma~\ref{lem-geometric-presentation} we obtain a new
surface $C_{e}$ of genus $e = d + 1$. The two open curves $r, s
\subset C_{d}\setminus D_{p}$ can be completed to smooth circles
 by adding arcs in
$C_{r,s}$ as shown on Figure~\ref{fig7}. Furthermore, we can enlarge the set of
relations $R$ to a new set of circles $S = \{ q \} \cup R \cup \{a, b
\}$ where $q \subset C_{r,s}$ is a circle isotopic to $r$ such that $r
\cap D_{p} = \varnothing$ and $a, b$ are the standard generators of
the fundamental group of $C_{p}$ (see Figure~\ref{fig7}).

\

\noindent
By construction we have a geometric presentation $\pi_{e} \to \Gamma$ 
whose kernel is generated as a normal subgroup by the elements of
$S$. By Lemma~\ref{lem-cyclic} the abelianization of the group
$\op{Map}_{S}$ is cyclic. In fact since we have ensured that $S$ is
big enough we have the following 

\begin{lem} \label{lem-cyclic-finite}
The abelianization of $\op{Map}_{S}$ is a finite cyclic group of order
dividing $10$.
\end{lem}
{\bf Proof.} The union
$(C_{r,s}\setminus D_{p})\cup C_{p}$ is a genus two handle on $C_{e}$
and therefore the inclusion $(C_{r,s}\setminus D_{p})\cup C_{p}
\subset C_{e}$ induces a natural homomorphism $h : \op{Map}_{2,1}
\to \op{Map}_{S}$. Moreover by construction the subgroup
$h(\op{Map}_{2,1}) \subset \op{Map}_{S}$ contains a non-trivial Dehn
twist - e.g. the twist along the circle $a \subset C_{e}$. Therefore
we conclude as in the proof of Lemma~\ref{lem-cyclic} that the
group $\op{Map}_{S}$ is generated by conjugates of the element 
$t_{a} \in h(\op{Map}_{2,1}) \subset \op{Map}_{S}$.

On the other hand it is known (see e.g. \cite{wajnryb}) that  the 
group $\op{Map}_{2,1}$ has a presentation with
several relations of degree zero and one relation of degree
ten. More precisely one has 
\[
\op{Map}_{2,1} = \left\langle t_{1}, \ldots , t_{5} \, \left| \,
\begin{minipage}[c]{3.5in} 
$t_{i}t_{j}t_{i} = t_{j}t_{i}t_{j}$ if $|i-j| = 1$, $t_{i}t_{j} =
t_{j}t_{i}$ if $|i-j| \neq 1$ and $(t_{1}t_{2}t_{3})^{4} = 
t_{5}(t_{4}t_{3}t_{2}t_{1}^{2}
t_{2}t_{3}t_{4})^{-1}t_{5}t_{4}t_{3}t_{2}t_{1}^{2}t_{2}t_{3}t_{4}$
\end{minipage}
\right. \right\rangle
\]
The braid relations  $t_{i}t_{j}t_{i} = t_{j}t_{i}t_{j}$ for $|i-j| =
1$ in this presentation show that any homomorphism from $\op{Map}_{2,1}$ to
an abelian group $G$ will have to map all $t_{i}$'s to the same
element $g \in G$. Moreover the  degree ten relation in the
presentation implies that $g^{10} = 1 \in G$ and so the abelianization
of $\op{Map}_{2,1}$ is isomorphic to ${\mathbb Z}/10$.

Due to this we see that the tenth power of any generator of
$\op{Map}_{S}$ will be conjugate to $t_{a}^{10}$ and so will belong to
the commutator subgroup $[\op{Map}_{S},\op{Map}_{S}]$. Thus by the
remark after the proof of Lemma~\ref{lem-cyclic} the abelianization of
$\op{Map}_{S}$ will be a cyclic group of order dividing ten. The lemma
is proven \hfill $\Box$

\medskip

Now we can finish the proof of the proposition. Let $m$ be the cardinality
of $S$ and let $s_{1}, \ldots, s_{m}$ be an ordering of the elements
in $S$. By Lemma~\ref{lem-cyclic-finite} we can find elements
$\xi_{i}, \eta_{i} \in \op{Map}_{S}$, $i = 1, \ldots, k$ such that 
the relation 
\begin{equation} \label{commutator}
\prod_{j = 1}^{m} t_{s_{j}}^{10} = \prod_{i=1}^{k} [\xi_{i},\eta_{i}]
\end{equation}
holds in $\op{Map}_{S}$. 

For every $i = 1, \ldots, n$ choose diffeomorphisms $\Xi_{i}, E_{i} \in
\op{Diff}^{+}(C_{e})$ representing the mapping classes  $\xi_{i},
\eta_{i}$. Let $K$ be a Riemann surface of genus $k$
with one boundary component. The fundamental group $\pi_{1}(K)$
is generated by $2k+1$ simple closed curves $\alpha_{1}, \beta_{1},
\ldots, \alpha_{k}, \beta_{k}, c$ satisfying the only relation 
\[
\prod_{i = 1}^{k}[\alpha_{i},\beta_{i}] = c
\]
and so we can construct a representation $\pi_{1}(K) \to
\op{Diff}^{+}(C_{e})$ given by $\alpha_{i} \mapsto \Xi_{i}$,
$\beta_{i} \mapsto E_{i}$, $c \mapsto \prod [\Xi_{i},E_{i}]$. Moreover
since we have the freedom of changing the diffeomorphisms $\Xi_{i}$
and $E_{i}$ up to isotopy we may assume without a loss of generality
that all the $\Xi_{i}$'s and $E_{i}$'s preserve a given symplectic
form on $C_{e}$. With this choice the group $\pi_{1}(K)$ acts
symplectically and freely on the product $\widetilde{K}\times C_{e}$
of the universal cover of $K$ and $C_{e}$ and by passing to the
quotient of this action we obtain a smooth symplectic fibration $Y \to
K$ with fiber $C_{e}$. The restriction of this fibration over the
boundary circle $c$ of $K$ is a principal $C_{e}$ fibration over
$S^{1}$ corresponding to the element $\prod_{i}[\xi_{i},\eta_{i}] \in
\op{Map}_{e}$. On the other hand exactly as in the proof of
Lemma~\ref{lem-tlf} we can construct a symplectic Lefschetz fibration
$u : U \to D$ whose restriction $U_{|\partial D}$ to the boundary
circle is a principal $C_{e}$ fibration over $S^{1}$ corresponding to
the element $\prod_{j} t_{j}^{10}$. Due to the relation
(\ref{commutator}) this implies that the restrictions of $Y$ and $U$
to $c \subset K$ and $\partial D \subset D$ respectively are
isomorphic $C_{e}$ bundles over a circle and so we can glue $Y$ and
$U$ along their boundaries to obtain a smooth orientable compact
fourfold $X$ which fibers over the genus $k$ curve $C_{k} = 
K\cup_{c=\partial D} D$ with a regular fiber of genus $e$. Since
the embeddings $U_{|\partial D}\subset U$ and $Y_{|c}\subset Y$ have
trivial normal bundles we can apply Weinstein's symplectic
neighborhood theorem \cite[Lemma~2.1]{GOMPF} 
to patch the symplectic forms on $U$ and
$Y$ over the gluing. Therefore the projection map $f : X \to C_{k}$ is
a symplectic Lefschetz fibration with $10m$ singular fibers and
vanishing cycles $s_{1}, \ldots, s_{m}$. In particular
$\ker[\pi_{1}(X) \to \pi_{1}(C_{k})]$ is isomorphic to the quotient of
$\pi_{e}$ by the normal subgroup generated by the $s_{j}$'s and so the
proposition is proven. \hfill $\Box$

\begin{rem} \label{rem-genus3}
It is clear from the proof of Proposition~\ref{prop-high} that with a
little extra care one can control effectively the genus of the fibers
of $f$ and the number of singular fibers of $f$ but not the genus of
the base $C_{k}$.

The trick of enlarging the geometric presentation by creating a genus
two handle seems to be really necessary since a priori one can 
guarantee only the existence of a non-trivial homomorphism 
from $\op{Map}_{1,1}$ to $\op{Map}_{d}$
which is not enough since the abelianization of $\op{Map}_{1,1}$ is
${\mathbb Z}$ as explained in the proof of
Corollary~\ref{identity}. On the other hand it is known that the group
$\op{Map}_{g,r}$ is perfect for $g \geq3$, $r \geq 0$
\cite{powell}. Hence same argument as in the proof of
Proposition~\ref{prop-high} shows that we can glue in an extra handle to
$C_{e}$  to obtain a non-trivial
homomorphism $\op{Map}_{3,1} \to \op{Map}_{d+2,1}$ and a geometric
presentation $\pi_{d+2} \to \Gamma$ whose kernel is generated as
a normal subgroup by a connected graph of circles $P$ with cardinality
$\# P = \# R + 5 = m + 2$. In this case already the Dehn twists $t_{s}$, $s
\in P$ themselves will be products of 
commutators in $\op{Map}_{P}$ and so we will
get a SLF with $m+2$ singular fibers.
\end{rem}

In fact by further enlarging the genus of the base one can modify the
fibration described in Remark~\ref{rem-genus3} to obtain 
a proof of Theorem~\ref{thm-high}. Actually we have the following precise
version of Theorem~\ref{thm-high}:

\begin{corr} \label{cor-one-fiber}
Let $\Gamma$ be a finitely presentable group with a given geometric
presentation $a: \pi_d \to \Gamma$. Then there exist a geometric
presentation $\pi_{d+2} \to \pi_{d} \stackrel{a}{\to} \Gamma$
and 
a symplectic Lefschetz fibration $f : Y \to C$ over a smooth compact
Riemann surface $C$
so that 
\begin{itemize}
\item the regular fiber of $f$ is diffeomorphic to $C_{d+2}$, 
\item $f$ has a unique singular fiber,
\item the
fundamental group of $X$ fits in a short exact sequence
\[
1 \to \Gamma \to \pi_{1}(X) \to \pi_{1}(C) \to 1
\]
\item the natural surjective map from the fundamental group of
the fiber of $f$ to $\Gamma$  coincides with the composition $\pi_{d+2}
\to \pi_{d} \stackrel{a}{\to} \Gamma$
\end{itemize}
\end{corr}
{\bf Proof.} We will use the notation of Remark~\ref{rem-genus3}.
As a first approximation to $f : X \to C$ we will use the genus $d+2$
SLF constructed in Remark~\ref{rem-genus3}. Observe that
since each $t_{s}$, $s \in P$ is a product of commutators in
$\op{Map}_{P}$ similarly to the proof of Theorem~\ref{thm-high}
we can replace a tubular neighborhood of the singular
fiber at which the circle $s$ vanishes with a smooth $C_{d+2}$
fibration over a high genus Riemann surface corresponding to a
representation of $t_{s}$ as a product of commutators. If we do this
for all elements in $S$ but one we will get a symplectic Lefschetz
fibration with one singular fiber and the same monodromy group as that
of the original fibration. Let $f : X \to C$ be this fibration and let
$q \in C$ be the only critical value of $f$. Let $r \in P$ be the
cycle vanishing at $q$. The map $f$ induces a
surjection on fundamental groups $f_{*} : \pi_{1}(X) \to \pi_{1}(C)$
and the kernel $\ker(f_{*})$ is naturally a quotient of the
fundamental group $\pi_{d+2}$ of the regular fiber of $f$. To
calculate $\ker(f_{*})$ note that we have a natural surjection
$\pi_{1}(X\setminus f^{-1}(q)) \to \pi_{1}(X)$ which restricts to the
quotient map $\pi_{d+2} \twoheadrightarrow \ker(f_{*})$. In
particular $\ker[\pi_{d+2} \twoheadrightarrow \ker(f_{*})]$
is precisely the subgroup in $\pi_{1}(X\setminus f^{-1}(q))$ which is
generated as normal subgroup (in $\pi_{1}(X\setminus f^{-1}(q))$) 
by the cycle $r \in \pi_{d+2} \subset \pi_{1}(X\setminus f^{-1}(q))$. 
On the other hand
the fact that $X\setminus f^{-1}(q) \to C\setminus\{q\}$ is a smooth
fibration identifies $\pi_{1}(X\setminus f^{-1}(q))$ with the
semi direct product $\pi_{1}(C\setminus\{q\})\ltimes_{\op{mon}}
\pi_{d+2}$ where $\op{mon} : \pi_{1}(C\setminus\{q\}) \to
\op{Aut}(\pi_{d+2})$ is the monodromy representation. By the
definition of a semi direct product the subgroup
$\pi_{1}(C\setminus \{q\}) \subset \pi_{1}(X\setminus f^{-1}(q))$
normalizes $\pi_{d+2} \subset \pi_{1}(X\setminus f^{-1}(q))$ and the
inner action of $\pi_{1}(C\setminus \{q\})$ on $\pi_{d+2}$ coincides
with the representation $\op{mon}$. Consequently $\ker[\pi_{d+2}
\twoheadrightarrow \ker(f_{*})]$ is the subgroup of $\pi_{d+2}$ which
is generated as a normal subgroup (in $\pi_{d+2}$) by the orbit of the
vanishing cycle $c \in \pi_{d+2}$ under the monodromy group of $f$.

The next step is to recall that by construction the monodromy group of
$f$ is exactly $\op{Map}_{P}$. Moreover since $P$ is a set of
relations in a geometric presentation we know that the Dehn twists
about any two circles in
$P$ are conjugate and hence for any $s \in P$ there exists an element
$\phi \in \op{Map}_{P}$ for which $t_{s} = \phi\circ t_{r} \circ
\phi^{-1} = t_{\phi(r)}$. However we have argued above that
$\phi(r) \in \ker(f_{*})$ for any $\phi \in \op{Map}_{P}$ and hence
$\ker{f_{*}}$ contains a circle which induces the same Dehn twist as
$s$. To finish the proof we only need to observe that for any surface 
$C_{g}$ and any circle
$c \subset C_{g}$ the Dehn twist $t_{c}$ is uniquely characterized
by the property that the maximal quotient of $\pi_{g}$ on which it 
induces the identity is the quotient of
$\pi_{g}$ by the normal subgroup generated by $c$. Thus if we have
two circles inducing the same Dehn twist the normal subgroups in $\pi_{g}$
generated by those circles coincide. In combination with the above
discussion this implies that $Q \subset \ker(f_{*})$ and so
$\ker(f_{*})$ coincides with the quotient of $\pi_{d+2}$ by the  
normal subgroup generated by the elements in $Q$. The corollary is
proven. \hfill $\Box$

\begin{rem} Observe that even though there are no obvious restrictions for
the Lefschetz fibrations constructed in Proposition~\ref{prop-high} to be
K\"{a}hler it will be very hard to construct projective examples like that. 
In general the algebraic Lefschetz fibrations have much bigger
monodromies and as a result the fundamental group of an algebraic
family $f : X \to C$ as above is an extension of  $\pi_{1}(C)$ with 
a rather small (in general trivial) group. The reason for this  
phenomenon is that as it follows from the above argument the only 
condition for the existence of a symplectic Lefschetz pencil with a 
given fundamental group is some equation in the 
mapping class group. On the other hand Hodge theory imposes many other
restrictions in the algebraic situation - for example it is shown in 
\cite{KPS} that if the first Betti number of $X$ is trivial, then 
the monodromy group of an algebraic Lefschetz pencil of big enough
degree cannot fix any loop on the fiber. An interesting question to
investigate is if this restriction  exists for SLP.
\end{rem}

\section{Examples}
\setcounter{figure}{0}

In this section we describe some  constructions of positive
relations among right Dehn twists that do not use the general method
employed in the proof of Theorem~\ref{main} but nevertheless lead
to SLF with interesting fundamental groups.

\subsection{Symplectic Lefschetz fibration with first Betti number
equal to one}

Our first construction produces a symplectic
Lefschetz fibration $f : X \to S^{2}$ with fibers of genus three having 
$b_{1}(X) = 1$. The first input for
the construction is provided by the following lemma.

\begin{lem} The group $\op{Map}_{2}$ contains five distinct right Dehn
twists $\{ t_{i} \}_{i=1}^{5} \subset \op{Map}_{2}$ satisfying the
positive relation
\[ 
(t_{1}^{2}\ldots t_{5}^{2})^{2} = 1.
\]
Moreover the $t_{i}$'s can be chosen to be twists along non-separating
circles in $C_{2}$. \label{lem-genus2relation}
\end{lem}
{\bf Proof.} We will construct the Dehn twists $t_{i}$ by exhibiting
an algebraic family of genus two curves which has a single 
non-reduced fiber and local monodromies given by $t_{i}^{2}$ at all
other singular fibers.

To achieve this  we start with the Hirzebruch surface ${\mathbb F}_{0} =
{\mathbb P}^{1} \times {\mathbb P}^{1}$ with the two standard
projections $p_{i} : {\mathbb F}_{0} \to {\mathbb P}^{1}$, $i =1,2$.
Consider ${\mathbb F}_{0}$ as a  ${\mathbb
P}^{1}$ bundle over ${\mathbb P}^{1}$ say via projecting onto the
first factor. Choose
six distinct points $q, q_{1}, q_{2}, \ldots, q_{5} \in {\mathbb
P}^{1}$ consider the divisors $F = p_{1}^{-1}(q)$, $H_{i} =
p_{2}^{-1}(q_{i})$, $i = 1, \ldots, 5$. If we denote by
$\Delta \subset {\mathbb F}_{0}$ the diagonal divisor in ${\mathbb
P}^{1} \times {\mathbb P}^{1}$, then the divisor $B := \sum_{i=1}^{5}H_{i}
+ F + \Delta$ belongs to the linear system $|{\mathcal O}_{{\mathbb
F}_{0}}(2,6)|$. Since the line bundle ${\mathcal O}_{{\mathbb
F}_{0}}(2,6)$ is divisible by two in $\op{Pic}({\mathbb F}_{0})$ the
divisor $B$ determines a two sheeted root cover $\alpha : S' \to {\mathbb
F}_{0}$ branched at $B$. The divisor $B$ has simple normal crossings
and so the covering $S'$ has an ordinary double
point at each ramification point of $\alpha$ sitting over a singularity of $B$.
Denote by $Q_{i} \in S'$ (respectively $P_{i} \in S'$), 
$i = 1, \ldots, 5$ the
ramification points of $\alpha$ sitting over the
intersection points of $\Delta$ (respectively $F$) with each of the
$H_{i}$'s and let $P_{0} \in S'$ be the ramification point of $\alpha$
projecting to the intersection point of $\Delta$ and $F$.

Let $\varepsilon : S \to S'$ be the minimal resolution of $S'$
obtained by blowing up the $P_{j}$'s and  the $Q_{i}$'s and let
$\mu_{j} = \varepsilon^{-1}(P_{j})$, $j = 0, \ldots, 5$ and $\nu_{i}
:= \varepsilon^{-1}(Q_{i})$, $i = 1, \ldots, 5$ be the corresponding
exceptional divisors. The composition $p := p_{1}\circ\alpha\circ
\varepsilon  : S
\to {\mathbb P}^{1}$ is a flat morphism having as critical values
exactly the points $q, q_{1}, \ldots, q_{5}$. For each $i = 1, \ldots,
5$ the fiber $p^{-1}(q_{i})$ is reduced and is the union of an
elliptic curve (the proper transform of $(p_{1}\circ
\alpha)^{-1}(q_{i}))$ and a $-2$ rational curve (the exceptional curve
$\nu_{i}$) meeting the elliptic curve transversally at two
points. Therefore the geometric vanishing cycle at each of the points
$q_{i}$ is a system of two isotopic non-intersecting circles
$c_{i}'\cup c_{i}'' \subset C_{2}$. This implies that the geometric 
monodromy transformation around each $q_{i}$ is given by the square of
the Dehn twist along $c_{i}'$ (or equivalently along $c_{i}''$).
Therefore in $\op{Map}_{2}$ we have a relation
\[
t_{1}^{2}\ldots t_{5}^{2} = \phi^{-1}
\] 
where $\phi \in \op{Map}_{2}$ is the local geometric monodromy around
the fiber $p^{-1}(q)$. Since $F = p_{1}^{-1}(q)$ is part of the branch
locus of the covering $\alpha$ the scheme-theoretic fiber $p^{-1}(q)$
will be nonreduced and will consist of one double rational component (whose
reduction maps one to one onto $F$) and six reduced rational
components (the exceptional divisors $\mu_{j}$). Finally observe that
if we pull back (and normalize)
the family $p : S \to {\mathbb P}^{1}$ by a double
cover $D := {\mathbb P}^{1} \to {\mathbb P}^{1}$ branched at $q$ and at one
extra point different from the $q_{i}$'s we will obtain a family over
$C$ with
a smooth total space whose fiber over the point $\tilde{q} \in C$ 
sitting over $q$ is reduced and has one component of genus 2 (the
double cover of $F$ branched at the points $P_{j}$, $j = 0,\ldots ,5$)
and six components (sitting over the $\mu_{i}$'s) which are smooth
rational curves of self-intersection $-1$. After contracting the six
$-1$ curves in the fiber we will obtain a genus two fibration which is
smooth at $q$ and hence has trivial monodromy at $q$. Thus the
monodromy transformation $\phi$ is of order two in $\op{Map}_{2}$ and
so the $t_{i}$'s satisfy
\[
(t_{1}^{2}\ldots t_{5}^{2})^{2} = 1.
\]
The lemma is proven. \hfill $\Box$

\begin{rem} (i) Clearly the construction used in the proof of
Lemma~\ref{lem-genus2relation} is very flexible and can 
be adapted to many different situations. 
One obvious generalization is the existence of $2n+1$ distinct Dehn
twists in $\op{Map}_{n}$ satisfying the relation
\[ 
(t_{1}^{2}\ldots t_{2n+1}^{2})^{2}=1.
\]  
\

\medskip

\noindent
(ii) Using the concrete nature of the construction in
Lemma~\ref{lem-genus2relation} one can obtain precise information
about the homology classes and the position 
of the vanishing cycles $c_{i}'$ and $c_{i}''$ (see
Lemma~\ref{lem-genus2character} below).  In fact one can show that 
the free loops $c_{i}'$, $i = 1, \ldots, 5$ generate the fundamental
group $\pi_{2}$ of the regular fiber of $S$ and that the Dehn twists
$t_{i}$, $i = 1, \ldots, 5$ generate the mapping class group
$\op{Map}_{2}$. 
\end{rem}

Before we state the second input to the construction we will need some
preliminaries. We will use the notation introduced in the proof of
Lemma~\ref{lem-genus2relation}. 
Fix a reference point $o \in {\mathbb P}^{1}\setminus \{q_{0},q_{1},
\ldots, q_{5}\}$ and let $S_{o} \cong C_{2}$ be the corresponding
fiber of $S$. Denote by $v_{i} = [c_{i}'] = [c_{i}''] \in
H_{1}(S_{o},{\mathbb Z})$, $i = 1, \ldots, 5$  the homology classes of
the vanishing cycles. Now we have

\begin{lem} \label{lem-genus2character} There exists a character 
$\chi : H_{1}(S_{o},{\mathbb Z}) \to {\mathbb Z}/2$ which is uniquely
determined by the property $\chi(v_{1}) = 0$ and 
$\chi(v_{i}) = 1$ for $i \neq 1$.
\end{lem}
{\bf Proof.} Since the rank of the free abelian group
$H_{1}(S_{o},{\mathbb Z})$ is four the existence of $\chi$ will be
proven if we can show that $v_{1},
v_{2}, v_{3}, v_{4}$ generate $H_{1}(S_{o},{\mathbb Z})$.

Note that by construction the family of hyperelliptic curves $p : S
\to {\mathbb P}^{1}$ is isotrivial and that the corresponding covering
of Weierstrass points splits into six components (the strict
transforms of the preimages of the divisors $H_{i}$ and the divisor
$\Delta$). Therefore in each pair $\{c_{i}', c_{i}'' \}$, $i = 1,
\ldots, 5$ the vanishing cycles $c_{i}'$ and $c_{i}''$ get
interchanged by the hyperelliptic involution and the image of the
curve $c_{i}'\cup c_{i}''$ via the projection $g := p_{2}\circ \alpha \circ
\varepsilon : S \to {\mathbb P}^{1}$ is a simple closed curve
$\gamma_{i} \subset {\mathbb P}^{1}\setminus \{q_{0},q_{1},
\ldots, q_{5}\}$ with the property that ${\mathbb P}^{1}\setminus
\gamma_{i}$ is a union of two open disks one of which contains
$\{q_{o}, q_{i} \}$ and the other contains $\{ q_{j} \}_{j \neq
0,i}$. 

\begin{figure}[!ht]
\begin{center}
\epsfig{file=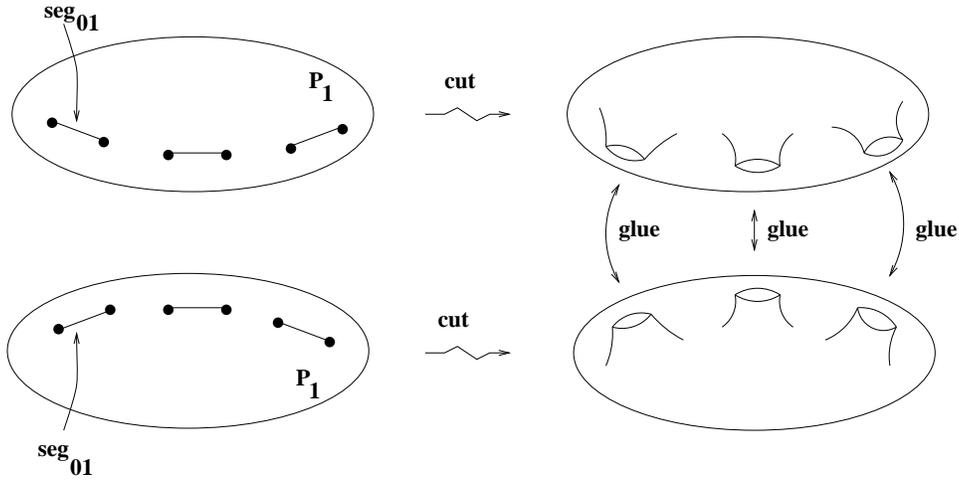,width=5in}
\end{center}
\caption{Gluing two ${\mathbb P}^{1}$'s into a hyperelliptic curve.}
\label{fig8} 
\end{figure}

Now using the fact that $g_{|S_{o}} : S_{o} \to {\mathbb P}^{1}$ is
the hyperelliptic covering we
can identify explicitly the homology classes $v_{i}$, $i = 1, \ldots,
5$ as follows. Choose three non-intersecting segments $\op{seg}_{01},
\op{seg}_{23}, \op{seg}_{45}$ in ${\mathbb
P}^{1}\setminus \{q_{0},q_{1},\ldots, q_{5}\}$ with endpoints $\{q_{0},
q_{1} \}$, $\{q_{2}, q_{3} \}$ and $\{q_{4}, q_{5} \}$
respectively. As usual the hyperelliptic involution allows one to
identify the genus two 
curve $S_{o}$ with the topological surface which is obtained by first
cutting ${\mathbb P}^{1}$ along the segments $\op{seg}_{01},
\op{seg}_{23}, \op{seg}_{45}$ and then gluing two copies of the
resulting surface along the opposite shores of the cuts (see
Figure~\ref{fig8}).

\

\begin{figure}[!ht]
\begin{center}
\epsfig{file=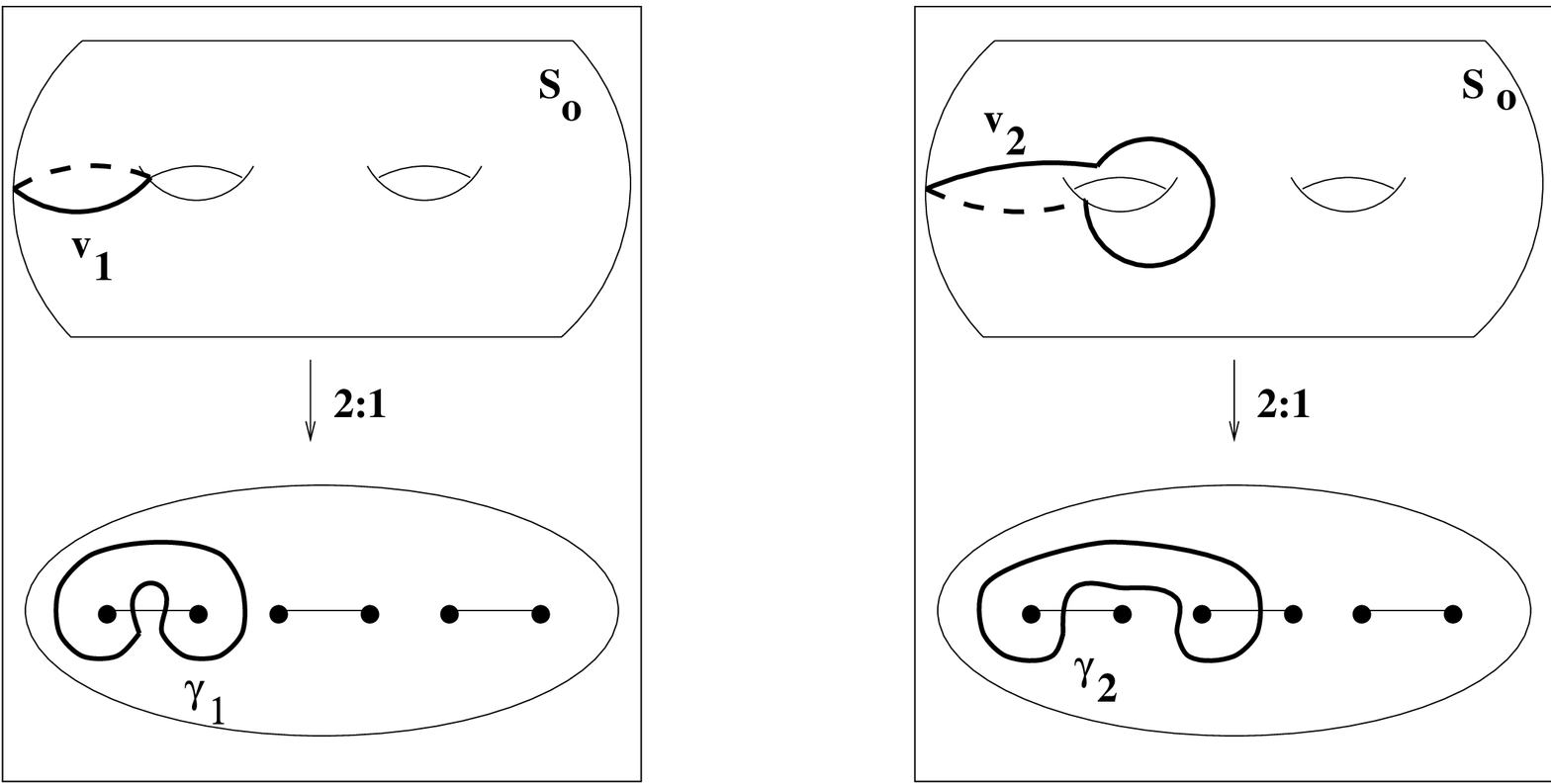,height=2.4in}

\smallskip

\epsfig{file=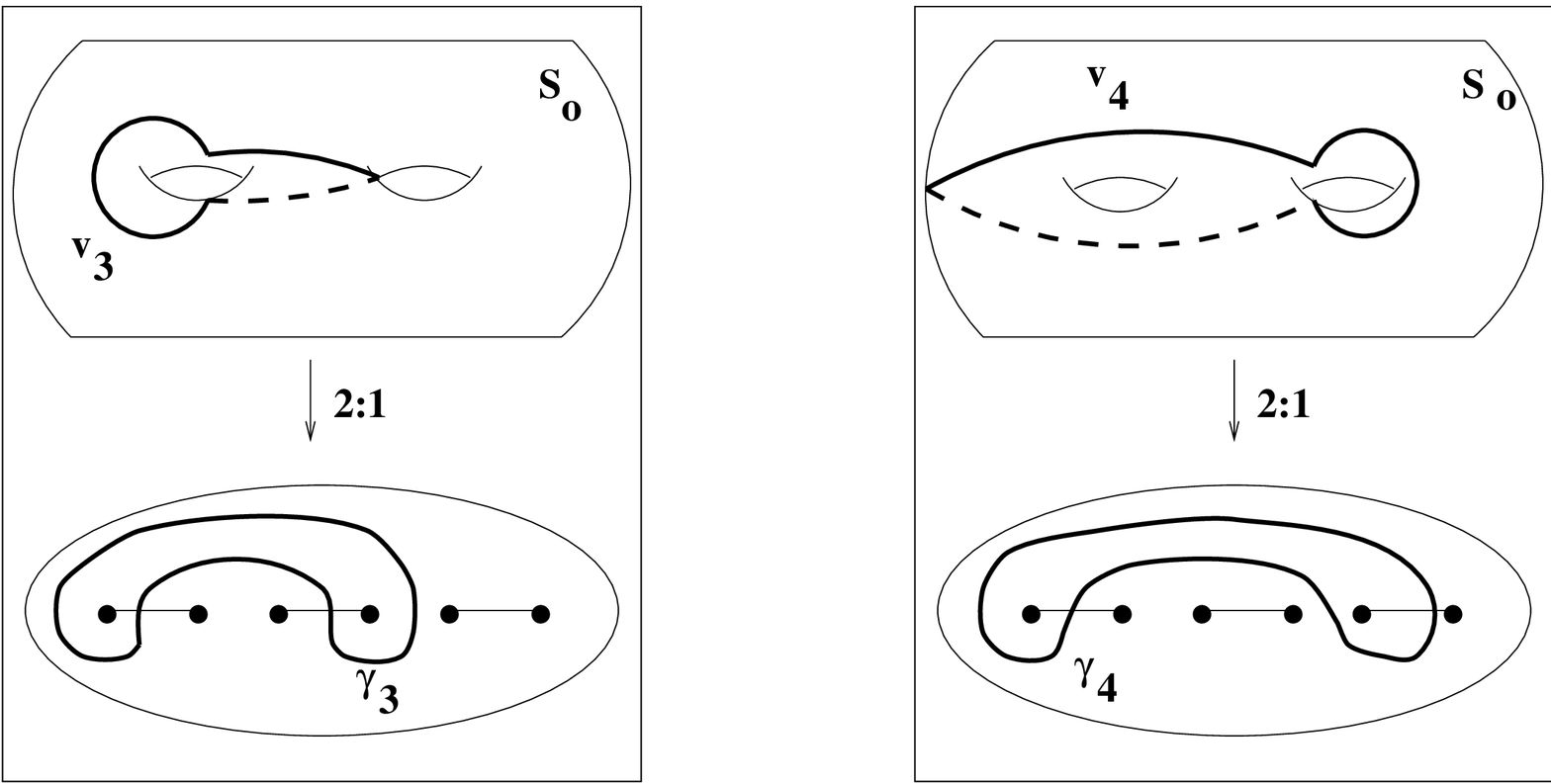,height=2.4in}

\smallskip

\epsfig{file=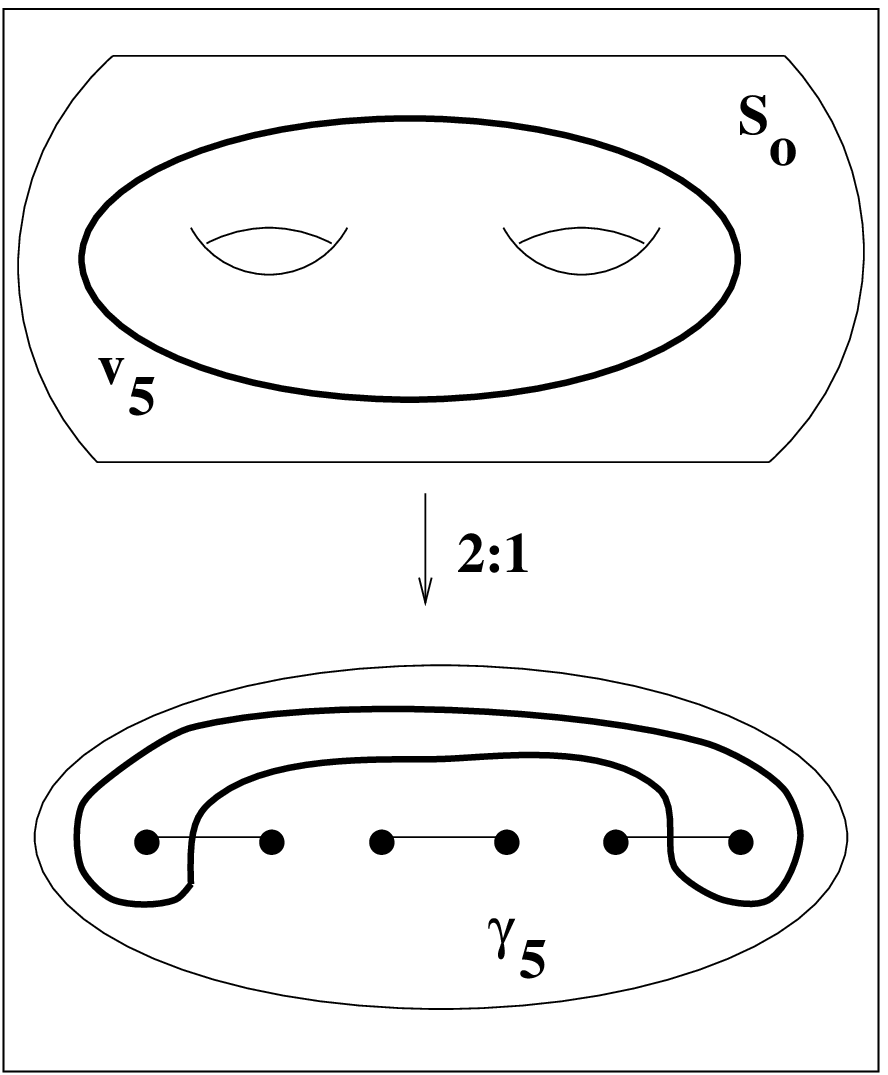,height=2.4in}
\end{center}
\caption{Loops representing the homology classes $\{ v_{i}
\}_{i=1}^{5}$.}  \label{fig9}
\end{figure}

\noindent
This explicit topological model of $S_{o}$ together with the above
characterization of the loops $\gamma_{i}$ shows that the classes
$v_{i}$, $i = 1, \ldots, 5$ are represented by the loops shown on
Figure~\ref{fig9}.

\

\begin{figure}[!ht]
\begin{center}
\epsfig{file=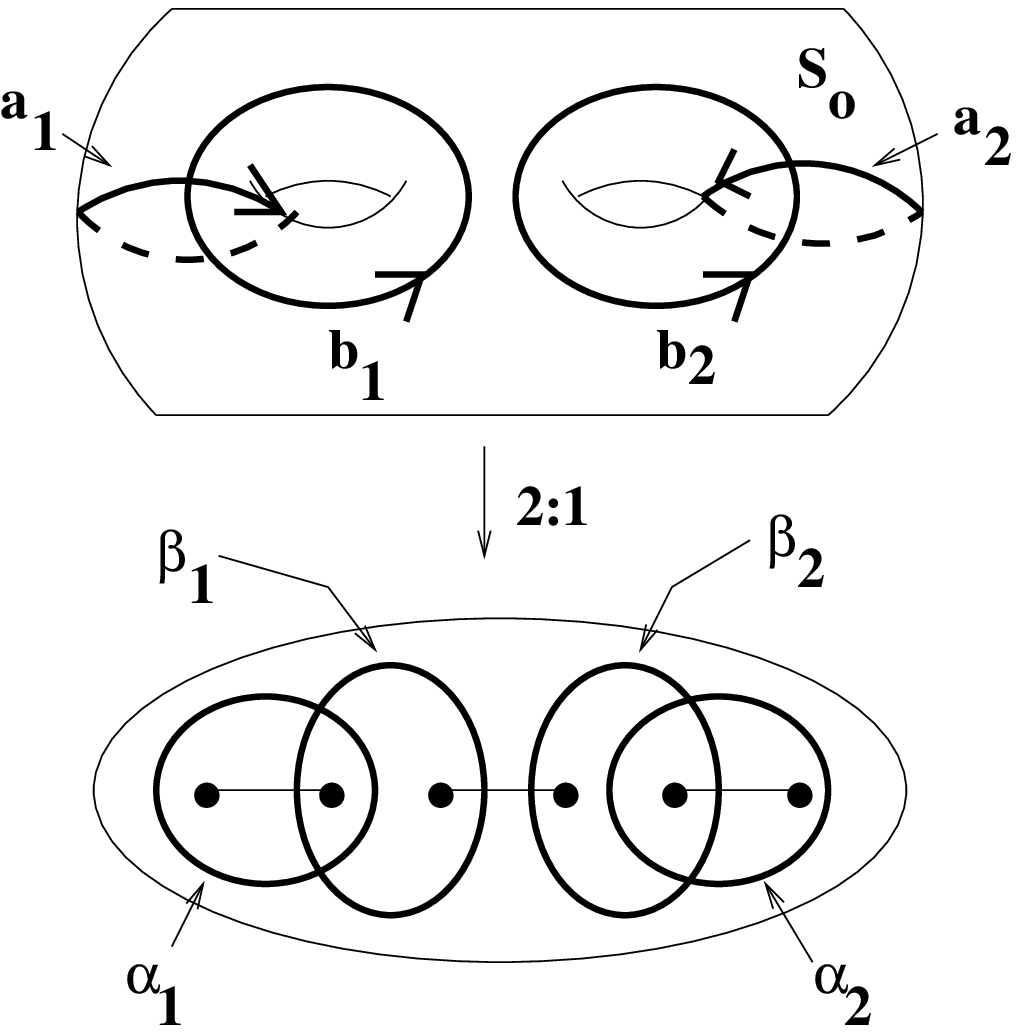,height=2.5in}
\end{center}
\caption{A standard basis in $H_{1}(S_{o}, {\mathbb Z})$.}  \label{fig10}
\end{figure}

\noindent
Let now $a_{1}, b_{1}, a_{2}, b_{2}$ be the standard basis of
$H_{1}(S_{o}, {\mathbb Z})$ which via the hyperelliptic map $g$ 
projects to the loops $\alpha_{1},
\beta_{1}, \alpha_{2}, \beta_{2}$ as depicted on Figure~\ref{fig10}. 

\

\noindent
Using the explicit realization of the $v_{i}$'s one calculates
$v_{1} = a_{1}$, $v_{2} = a_{1} - b_{1}$, $v_{3} = -a_{1} - b_{1} + a_{2}$,
$v_{4} = - b_{1} - a_{2} + b_{2}$, $v_{5} = b_{1} + b_{2}$.
For this calculation we have used the standard orientations of the
basis elements shown on Figure~\ref{fig10} and suitable orientations of the
$v_{i}$'s. Note that since we are ultimately interested in the
homology of $S_{o}$ only modulo 2 the ambiguity introduced by the
choice of orientations of the individual $v_{i}$'s will not affect the
anything.

The above formulas show that  the desired character $\chi$ is
just the character for which $a_{1} \mapsto 1, b_{1} \mapsto 0, a_{2}
\mapsto 1, b_{2} \mapsto 1$. The lemma is proven. \hfill $\Box$

\bigskip

Now we move to the  construction of a family of genus 3 curves $X$ with 
$b_{1}(X) = 1$.

\medskip

\noindent
Let $C_{2}$ be a smooth surface of genus two and let
$c_{1}', \ldots, c_{5}' \subset C_{2}$ be circles giving the Dehn twists
$t_{1}, \ldots, t_{5}$ constructed in Lemma~\ref{lem-genus2relation}. 
Let $d : D \to C_{2}$ be the
unique unramified double cover of $C_{2}$ corresponding to the
character $\chi$ from Lemma~\ref{lem-genus2character}. By the
definition of $\chi$ we know that $d^{-1}(c_{i}') = w_{i}$
is connected for $i \neq 1$ and $d^{-1}(c_{1}) = u\cup v$ is a union of two
non-intersecting circles which get interchanged by the involution on $D$.
Therefore in $\op{Map}_{3}$ we get a relation
\[ 
(t_{u}t_{v}t_{w_{2}}t_{w_{3}}t_{w_{4}}t_{w_{5}})^{2} = 1
\]
Let $f : X \to S^{2}$ be the symplectic Lefschetz fibration determined
by this relation as in Lemma~\ref{lem-tlf}.

The homology classes of the $w_{i}$ clearly generate the rank four
subspace of $H_{1}(D,{\mathbb Q})$ consisting of cycles invariant
under the covering involution for the cover $d: D \to
C_{2}$. Furthermore by construction the cycle $u$ is not in the span
of the $w_{i}$'s but $u+v$ is invariant under the involution and so
$w_{1}, w_{2}, w_{4}, w_{5}, u, v$ span a five dimensional subspace in
$H_{1}(D,{\mathbb Q})$. Since by Lemma~\ref{lem-tlf} (c) we know that 
$H_{1}(X,{\mathbb Q})$ is the quotient of $H_{1}(D,{\mathbb Q})$
by the linear span of the vanishing cycles we conclude that $b_{1}(X)
= 1$ as desired.

\begin{rem} The family $f : X \to S^{2}$ just constructed  is an
example of a symplectic Lefschetz fibration whose monodromy group does
not act semi simply on the first cohomology of the fiber. This can be
seen directly from the explicit description of the Dehn twists
$t_{w_{i}}$, $t_{u}$ and $t_{v}$ given above. In fact it is proven in
\cite{KPS} that as long as the total space of a SLF has an odd first
Betti number the monodromy action on the first cohomology of the fiber
has a non-trivial unipotent radical.
\end{rem}

\begin{rem} Explicit examples of symplectic Lefschetz fibrations with
prescribed first homology were also constructed by Ivan Smith
\cite{smith}. In particular he shows that any abelian group which is a
quotient of a free group on $g$ generators can be realized as the
first homology of a SLF with fiber genus $2g$.
\end{rem}

\subsection{Symplectic
Lefschetz fibration with  fundamental group equal to
${\mathbb Z}$.}

In this section we will use a slight modification of the proof of
Proposition~\ref{prop-main} to construct an example of a  symplectic
Lefschetz fibration with  fundamental group equal to
${\mathbb Z}$.

Let $E$ be an elliptic curve. Consider the
complex projective surface $V := E\times {\mathbb P}^1$ and let 
$p_{E}: V \to E$ be the natural projection. Let $D_{\nu} \subset V$ 
be the section of $p_{E}$
given by the graph of the natural degree two covering $\nu : E \to {\mathbb
P}^1$. Fix a point $p \in
{\mathbb P}^{1}$ which is not a branch point for the covering 
$\nu : E \to {\mathbb
P}^1$ and let $D_{p} = E \times \{ p \}$ be the corresponding
section of $p_{E}$.

\begin{figure}[!ht]
\begin{center}
\epsfig{file=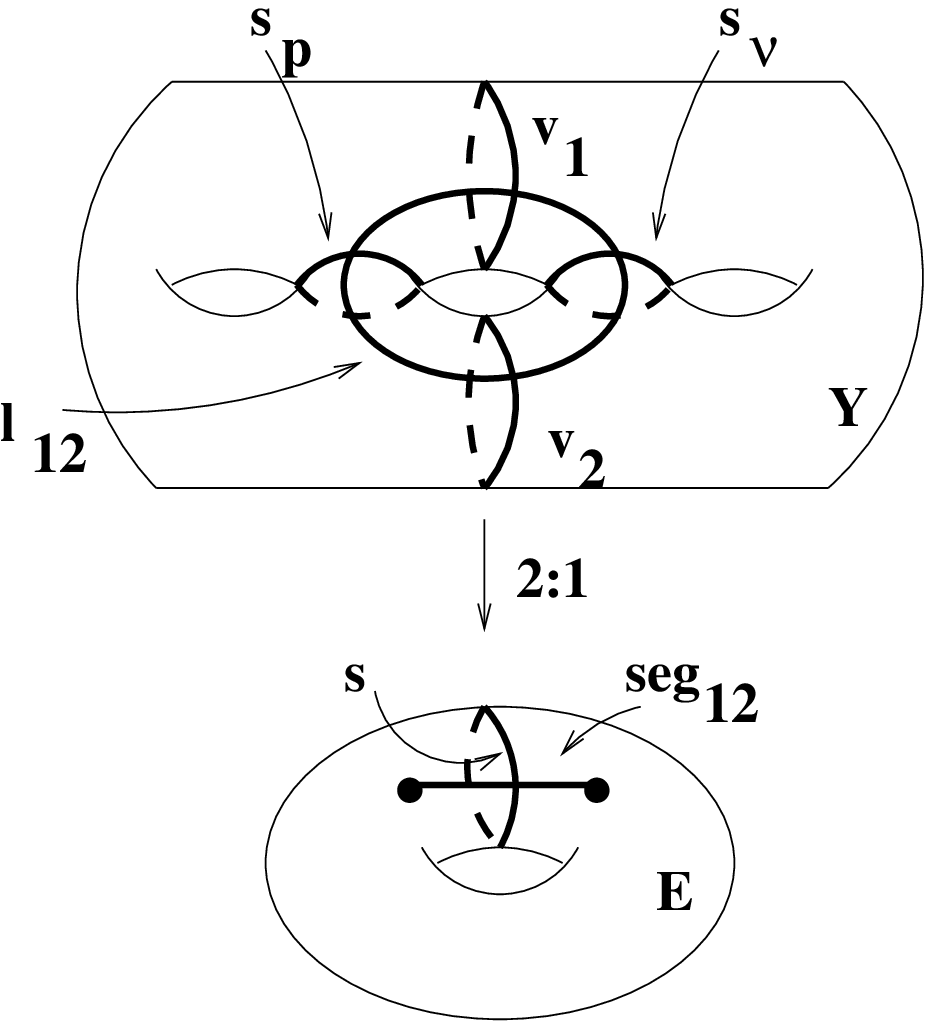,height=2.5in}
\end{center}
\caption{The covering $p_{E|Y} : Y \to E$.}  \label{fig11}
\end{figure}

Consider a divisor $D = D_{p} + D_{\nu}$.
By construction $D$ is a strict normal crossing curve with exactly two
singular points $\{p_{1}, p_{2}\} = D_{p}\cap D_{\nu}$.
Suppose that the linear system $|D|$ contains a smooth divisor $Y$.
We will assume that $Y$ is close enough to $D$ in $|D|$ so that there
is a well defined deformation retraction $cr : Y \to D$.
The curve $Y$ is of genus $3$ and the natural map $p_{E|Y} : Y \to E$
is a two sheeted covering branched at four points which get glued in
pairs to the points
$p_{E}(p_{1})$ and $p_{E}(p_{2})$ respectively when we deform $Y$ to
$D$.  Let $\op{seg}_{12}$ be a contractible segment in $E$
connecting $p_{E}(p_{1})$ and $p_{E}(p_{2})$ and let $s \subset E$ be
a non separating circle intersecting $\op{seg}_{12}$ at a single point
(see Figure~\ref{fig11}).
Denote by $S$ the set of geometric vanishing cycles for an algebraic
pencil ${\mathbb P}^{1} \subset |D|$ containing both $Y$ and $D$.
As in the proof of Proposition~\ref{prop-main} consider the system of
circles in minimal position
\[
L := \{ l_{12}, s_{p}, s_{\nu} \} \cup S \subset Y
\]
where $l_{12} := p_{E|Y}^{-1}(\op{seg}_{12})$ is the double cover 
of $\op{seg}_{12}$ and $s_{p}$ and $s_{\nu}$ are the
preimages of $s$ in $E_{p}$ and $E_{\nu}$ respectively, viewed as cycles
on $Y$ via $cr$ (see Figure~\ref{fig11}). By definition the set $S$ 
contains two cycles 
$v_{1}$ and $v_{2}$ corresponding to
$p_{1}$ and $p_{2}$ respectively 
(i.e. the cycles contracted by $cr$) and so $L$ is graph connected.

\

\noindent
The same argument as in the proof of Proposition~\ref{prop-main} shows
that the quotient $\Gamma$ of $\pi_{1}(Y)$ by the normal subgroup generated by
the circles in $L$ is isomorphic to the quotient of $\pi_{1}(E)$ by
the normal subgroup generated by $s$, i.e. is isomorphic to ${\mathbb
Z}$. Furthermore since $L$ is graph connected the combination of
Lemma~\ref{lem-tlf} and the ``only if'' part of
Lemma~\ref{lem-negative} a SLF with
monodromy group $\op{Map}_{L}$, regular fiber isomorphic to $Y$ and
fundamental group $\Gamma$.

Therefore in order to finish the construction we only need to show
that the linear system $|D|$ contains a smooth divisor. Let ${\mathbb F}_{0} :=
{\mathbb P}^{1}\times {\mathbb P}^{1}$ and let $a := \nu\times \op{id}
: V \to {\mathbb F}_{0}$ be the double cover induced from $\nu$. As usually
${\mathcal O}_{{\mathbb F}_{0}}(m,n)$ will denote the line bundle 
$p_{1}^{*}{\mathcal O}_{{\mathbb P}^{1}}(m)\otimes p_{2}^{*}{\mathcal
O}_{{\mathbb P}^{1}}(n)$ where $p_{i} : {\mathbb F}_{0} \to {\mathbb
P}^{1}$, $i = 1,2$  are the projections on the first and second
factors respectively. By definition we have $D_{p} \in
|a^{*}{\mathcal O}_{{\mathbb F}_{0}}(0,1)|$ and $D_{\nu} \in
|a^{*}{\mathcal O}_{{\mathbb F}_{0}}(1,1)|$. Thus $D \in
|a^{*}{\mathcal O}_{{\mathbb F}_{0}}(1,2)|$ and so we need to find a
smooth divisor in the linear system $|a^{*}{\mathcal O}_{{\mathbb
F}_{0}}(1,2)|$. But the double cover $a : V \to {\mathbb F}_{0}$ is a
root cover branched along a section of ${\mathbb O}_{{\mathbb
F}_{0}}(4,0)$ and hence $a_{*}{\mathcal O}_{V} = {\mathcal O}\oplus
{\mathcal O}(-2,0)$. This yields
\begin{align*}
H^{0}(V,a^{*}{\mathcal O}_{{\mathbb F}_{0}}(1,2)) & = 
H^{0}({\mathbb F}_{0},a_{*}a^{*}{\mathcal O}_{{\mathbb F}_{0}}(1,2))
\\
& = H^{0}({\mathbb F}_{0},{\mathcal O}_{{\mathbb F}_{0}}(1,2)\otimes
({\mathcal O}_{{\mathbb F}_{0}}\oplus 
{\mathcal O}_{{\mathbb F}_{0}}(-2,0))) \\
& = H^{0}({\mathbb F}_{0},{\mathcal O}_{{\mathbb F}_{0}}(1,2))
\end{align*}
and so every member of the linear system $|D|$ is a pullback of a
divisor in $|{\mathcal O}_{{\mathbb F}_{0}}(1,2)|$. Finally since
${\mathcal O}_{{\mathbb F}_{0}}(1,2)$ is obviously base point free the
general point in the five dimensional projective space $|{\mathcal
O}_{{\mathbb F}_{0}}(1,2)|$ will represent a smooth rational curve on
${\mathbb F}_{0}$ which is a double cover of ${\mathbb P}^{1}$ via the
projection $p_{1}$ and which intersects the branch divisor of $a$ in
eight distinct points. Therefore the preimage of such a curve in $V$
is a smooth curve $Y$ of genus three which finishes the construction.

\begin{rem}
As an alternative to the previous construction one can use our main
theorem in the following way. Start with a finitely presentable group
a $\Gamma$ and let 
$a_{1} : \pi_g \twoheadrightarrow \Gamma$ and $a_{2}:\pi_g
\twoheadrightarrow \Gamma$ be two presentations. Use
Theorem~\ref{main} to construct two SLF $f_{1} : X_{1}
\to S^{2}$ and $f_{2} : X_{2} \to S^{2}$ of the same genus $h > g$ 
so that $\pi_{1}(X_{i}) = \Gamma$ and the natural maps from the
fundamental groups of the general fibers of $f_{i}$ onto $\Gamma$
factor as $\pi_{h} \to \pi_{g} \stackrel{a_{i}}{\to} \Gamma$ for $i =
1,2$. If we now glue the two pencils $X_{1}$ and $X_{2}$ along a
regular fiber we obtain a new SLF $f : X \to S^{2}$ of genus $h$ which
by van Kampen's theorem has fundamental group
$\pi_{g}/(\ker(a_{1})\ast \ker(a_{2}))$.
In particular by gluing two genus two pencils with 
fundamental groups ${\mathbb Z}^{2}$ one can get a genus two Lefschetz
pencil with fundamental group  ${\mathbb Z}$. In fact this procedure
is the baby version of a sophisticated technique employed by Ivan
Smith \cite{smithgenus2} who showed that every quotient of ${\mathbb
Z}^{\oplus 2}$ can be realized as the fundamental group of a
symplectic  Lefschetz fibration.

Similar  series of very elegant examples of genus two SLF with 
fundamental groups  ${\mathbb Z}\oplus {\mathbb Z}/n$ was constructed
in \cite{ozbagci-stipsicz}.
\end{rem}

\begin{rem} It was pointed out to us by Ron Stern that similarly to 
\cite{ozbagci-stipsicz} and \cite{smithgenus2}
the  examples with a fundamental group ${\mathbb Z}$ and $b_{1} = 1$
we construct above are all not even homotopic to complex surfaces. Indeed 
from the classification of complex surfaces it follows that if the
above surfaces are complex they are either secondary Kodaira surfaces,
class VII surfaces or elliptic surfaces. The fact that our examples
are symplectic excludes the first two possibilities. The third
possibility is ruled out by the observation that if this is an elliptic
fibration then it should be over $S^{2}$. But then $b_{1} \neq 1$
as it follows from \cite[Section~2.2]{FM}.
\end{rem}

\section{Sections in symplectic Lefschetz fibrations} \label{spheres}

In this section we study the relation between  numerical properties of the 
sections in a symplectic Lefschetz fibration $f : X \to S^{2}$ 
and the geometric monodromy of $f$. 

\subsection{Mapping class groups in genus one} \label{subsec-mapone}

First we recall some well known
facts on the structure of the mapping class groups in genus one (see
e.g. \cite{BIR,gervais}).

Since any translation on a torus is homotopic to the identity, the 
natural map $\op{Map}_{1}^{1} \to \op{Map}_{1}$ is an
isomorphism. Furthermore $\op{Map}_{1}^{1} = SL(2,{\mathbb Z})$ and
coincides with the group of linear automorphisms of the two
dimensional torus. The group $\op{Map}_{1,1}$ is naturally identified
with the mapping class group of a two dimensional torus with one
hole. As in Remark~\ref{remark-extension} the natural
forgetful map $\op{Map}_{1,1} \to \op{Map}^{1}_{1}$ realizes
$\op{Map}_{1,1}$ as a central extension
\begin{equation} \label{metaplectic}
0 \to {\mathbb Z} \to \op{Map}_{1,1} \to SL(2,{\mathbb Z}) \to 1,
\end{equation}
where the kernel ${\mathbb Z}$ is generated by the right twist $t_{c}$ along
the boundary circle $c$ of the hole.

On the other hand the group $\op{Map}_{1,1}$ admits a presentation
\begin{equation} \label{presentation}
\op{Map}_{1,1} = \langle t_{a}, t_{b} | t_{a}t_{b}t_{a} =
t_{b}t_{a}t_{b} 
\rangle,
\end{equation}
with generators $t_{a}$ and $t_{b}$ corresponding to the Dehn twist along a
standard symplectic basis of cycles $H_{1}(C_{1},{\mathbb Z}) = {\mathbb
Z}a\oplus {\mathbb Z}b)$ on the non-punctured torus $C_{1}$. In particular
under the natural map  $\op{Map}_{1,1} \to SL(2,{\mathbb Z})$ we have
\[
t_{a} \mapsto \begin{pmatrix}
1 & 1 \\ 0 & 1 \end{pmatrix}, t_{b} \mapsto 
\begin{pmatrix} 1 & 0 \\ -1 & 1
\end{pmatrix}.
\] 
From this is clear that the element
$(t_{a}t_{b})^{3}$ is central and maps to $- 1 \in SL(2,{\mathbb Z})$. 
Thus $(t_{a}t_{b})^{6}$ generates  ${\mathbb Z} = \ker[\op{Map}_{1,1} 
\to SL(2,{\mathbb Z})]$, i.e. $(t_{a}t_{b})^{6} = t_{c}$ (see also
\cite[Theorem~4.3]{ivanov-mccarthy}).  

The central extension (\ref{metaplectic}) corresponds to an element 
$\tau \in H^{2}(SL(2,{\mathbb Z}),{\mathbb Z})$. Since $SL(2,{\mathbb
Z})$ can be identified with the fundamental group of the moduli stack
${\mathcal M}_{1}^{1}$ of elliptic curves, the element $\tau$ can be 
interpreted as the first Chern class of a line bundle on ${\mathcal
M}_{1}^{1}$. In fact Mumford had shown
\cite[Main~Theorem]{mumford-picard} that
$\tau$ generates
$\op{Pic}({\mathcal M}_{1}^{1}) \cong {\mathbb
Z}/12$ and that $\tau$ corresponds to the line bundle on
${\mathcal M}_{1}^{1}$ which to every flat family $p: E \to S$ of elliptic
curves associates the element $\bar{\omega}_{f} \in \op{Pic}(S)$ where
$\bar{\omega}_{f}$ is the pullback of the relative dualizing sheaf of
$f$ via the zero section.

More algebraically the central extension (\ref{metaplectic}) can be
described as follows. Consider the universal covering
$\widetilde{SL}(2,{\mathbb R}) \to SL(2,{\mathbb
R})$. Since $\pi_{1}(SL(2,{\mathbb R})) \cong {\mathbb Z}$ is a
central subgroup in $\widetilde{SL}(2,{\mathbb R})$ we can take the
preimage $\widetilde{SL}(2,{\mathbb Z})$ of $SL(2,{\mathbb Z})$ in 
$\widetilde{SL}(2,{\mathbb R})$. By construction there is a natural
central extension
\begin{equation} \label{meta}
0 \to {\mathbb Z} \to \widetilde{SL}(2,{\mathbb Z}) \to SL(2,{\mathbb
Z}) \to 1
\end{equation}

One has the following simple but somewhat tedious lemma:

\begin{lem} \begin{list}{{\em (\alph{inner})}}{\usecounter{inner}}
\item The group $\op{Map}_{1,1}$ is isomorphic to 
$\widetilde{SL}(2,{\mathbb Z})$. 
\item The elements $t_{a}, t_{b} \in \op{Map}_{1,1}$ are conjugate in 
$\op{Map}_{1,1}$.
\end{list} \label{lem-sl2}
\end{lem}

\noindent
As an immediate corollary we get 

\begin{corr} The subsemigroup of $\op{Map}_{1,1}$ generated by the
conjugates of the element $t_{a}$ does not contain the identity
element.
\label{identity}
\end{corr}

The proofs of these statements will be given in
Sections~\ref{subsec-univ-covers} and \ref{subsec-sections} respectively.
The proof of Lemma~\ref{lem-sl2} is based on the
standard description \cite[Section~1.8-1.9]{lion-vergne}
of the universal cover of $SL(2,{\mathbb R})$
which we recall next.

\subsection{Universal covers of symplectic groups} \label{subsec-univ-covers}

Let $V$ be a two dimensional real vector space with basis $p, q \in V$
and coordinate functions $x, y \in V^{\vee}$. The group $SL(2,{\mathbb
R})$ is naturally identified with the group of linear automorphisms of
$V$ preserving the standard symplectic form $w := x\wedge y \in
\wedge^{2}V^{\vee}$. Let $\Lambda$ denote the Grassmanian of all
Lagrangian subspaces in $V$. Since $\dim_{{\mathbb R}}(V) = 2$ we have
$\Lambda = Gr(1,V) \cong {\mathbb R}{\mathbb P}^{1} \cong S^{1}$. In
fact any line $\ell \subset V$ is uniquely determined by the angle
$(\theta \mod \pi)$ where $\theta$ is the angle $\ell$ forms with the
$x$-axis. If we identify $V \cong {\mathbb C}$ via $(x,y) \mapsto x+
iy$, then the identification $u : \Lambda \widetilde{\to} S^{1}$ is
given explicitly as $u(\ell) = e^{2i\theta}$, where $\ell = {\mathbb
R}e^{i\theta}p$. The symplectic group $SL(2,{\mathbb R})$ acts on
$\Lambda$ and this action lifts to a well defined continuous action of
$\widetilde{SL}(2,{\mathbb R})$ on the universal cover
$\widetilde{\Lambda} \cong {\mathbb R}$ of $\Lambda$ which essentially
determines $\widetilde{SL}(2,{\mathbb R})$

Explicitly the group $\widetilde{SL}(2,{\mathbb R})$ can be described
in terms of the $SL(2,{\mathbb R})$ action on $\Lambda$ by means of
the Maslov index \cite{capell, lion-vergne}. 
Recall \cite[Section~1.5]{lion-vergne} that for a
real symplectic vector space $(V,w)$ with a Lagrangian Grassmanian
$\Lambda$, the corresponding Maslov index is the function $\tau :
\Lambda^{3} \to {\mathbb Z}$ defined as follows. For each triple 
$(\ell_{1},\ell_{2},\ell_{3})$ of Lagrangian subspaces of $V$ consider
the quadratic form $Q_{\ell_{1}\ell_{2}\ell_{3}}(x)$ on the vector
space $\ell_{1}\oplus\ell_{2}\oplus\ell_{3}$ given by
\[
Q_{\ell_{1}\ell_{2}\ell_{3}}(x_{1}\oplus x_{2}\oplus x_{3}) =
w(x_{1},x_{2}) + w(x_{2},x_{3}) + w(x_{3},x_{1}).
\]
Next put $\tau(\ell_{1},\ell_{2},\ell_{3}) := \text{ the signature of
} Q_{\ell_{1}\ell_{2}\ell_{3}}$. By construction the function $\tau$
is antisymmetric in the three arguments and is invariant under the
diagonal action of $Sp(V,w)$ on $\Lambda^{3}$. Moreover if $\ell \in
\Lambda$ is a fixed Lagrangian subspace of $V$, then the ${\mathbb
Z}$-valued function
\[
\xymatrix@R=1pt{
\tau_{\ell} : & Sp(V,w)\times Sp(V,w) \ar[r] & {\mathbb Z} \\
& (g,h) \ar[r] & \tau(\ell,g\ell,gh\ell),
}
\]
satisfies (see \cite[Lemma~1.6.13]{lion-vergne})
\[
\tau_{\ell}(g_{1}g_{2},g_{3}) + \tau_{\ell}(g_{1},g_{2}) = 
\tau_{\ell}(g_{1},g_{2}g_{3}) + \tau_{\ell}(g_{2},g_{3}).
\]
Thus $\tau_{\ell}$ is a cocycle for the group $Sp(V,w)$ with coefficients in
the trivial module ${\mathbb Z}$ and so determines a central extension
\begin{equation} \label{gell}
0 \to {\mathbb Z} \to G_{\ell} \to Sp(V,w) \to 1.
\end{equation}
As a set $G_{\ell} = Sp(V,w)\times {\mathbb Z}$ with the group
structure being given by 
\begin{equation} \label{grouplaw}
(g_{1},n_{1})\cdot(g_{2},n_{2}) = (g_{1}g_{2}, n_{1} + n_{2} +  
\tau_{\ell}(g_{1},g_{2})).
\end{equation}
The group $G_{\ell}$ acts naturally on the universal cover
$\widetilde{\Lambda}$ of $\Lambda$. It turns out
\cite[Section~1.9]{lion-vergne} that there is a unique topology on
$G_{\ell}$ for which this action becomes continuous, and so $G_{\ell}$
has a natural structure of a Lie group. Furthermore the universal
covering group $\widetilde{SP}(V,w)$ of $Sp(V,w)$ is naturally
identified with the identity component of $G_{\ell}$. In fact there is
a canonical character $s : G_{\ell} \to S^{1}$ of order four so that 
$\widetilde{SP}(V,w) = \ker(s)$.

This description of $\widetilde{SP}(V,w)$ can be made completely
explicit for the case of a two dimensional symplectic vector space $(V,w)$
considered above. Indeed the Maslov index map $\tau$ admits a very
concrete description in this case. The natural (counterclockwise)
orientation of $S^{1} \subset {\mathbb R}^{2}$ gives a cyclic ordering
for any triple of points in $S^{1}$. Let $\ell_{1},\ell_{2},\ell_{3}
\in \Lambda$ be three lines in $V$. Then the Masloiv index
$\tau(\ell_{1},\ell_{2},\ell_{3})$ is given by
\[
\tau(\ell_{1},\ell_{2},\ell_{3}) = \left\{ 
\begin{array}{rl}
0 & \text{ if } \ell_{1},\ell_{2},\ell_{3} \text{ are not all distinct
} \\
1 & \text{ if } u(\ell_{2}) \text{ is between } u(\ell_{1}) \text{ and
} u(\ell_{3})\\
-1 & \text{ if } u(\ell_{2}) \text{ is between } u(\ell_{3}) \text{ and
} u(\ell_{1})
\end{array}
\right.
\]
Geometrically this just means that $\tau(\ell_{1},\ell_{2},\ell_{3}) =
1$ if $\ell_{2}$ is inside the angle $(\mod \pi)$ formed by $\ell_{1}$
and $\ell_{3}$ and $\tau(\ell_{1},\ell_{2},\ell_{3}) = -1$ if
$\ell_{2}$ is outside that angle.

Choose $\ell := {\mathbb R}\, p$ and consider the central extension
(\ref{gell}). The function $s : G_{\ell} \to S^{1}$ defined by
\[
s(\begin{pmatrix} a & b \\ c & d \end{pmatrix}, n) := \left\{
\begin{array}{ll}
\op{sgn}(c)i^{n+1} & \text{ if } c \neq 0 \\
\op{sgn}(a)i^{n} & \text{ if } c = 0 
\end{array}
\right.
\] 
is a character of $G_{\ell}$ and $\widetilde{SL}(2,{\mathbb R}) =
\ker(s)$. In other words we have
\begin{equation} \label{usl2}
\widetilde{SL}(2,{\mathbb R}) = \left\{
(\begin{pmatrix} a & b \\ c & d \end{pmatrix}, n) \in SL(2,{\mathbb
R})\times {\mathbb Z} 
\left|
\begin{minipage}[c]{3in}
If $c = 0$, then $n$ is even and $\op{sgn}(a) = (-1)^{n/2}$ and if $c
\neq 0$, then $n$ is odd and $\op{sgn}(c) = (-1)^{(n+1)/2}$
\end{minipage}
\right.
\right\}
\end{equation}
with a group law given by the formula (\ref{grouplaw}).

\bigskip

Using this explicit description we can prove Lemma~\ref{lem-sl2}

\bigskip

\noindent
{\bf Proof of Lemma~\ref{lem-sl2}.} Put
\[
A := \begin{pmatrix} 1 & 1 \\
0 & 1 \end{pmatrix} \;\; B := \begin{pmatrix} 1 & 0 \\ -1 & 1
\end{pmatrix}, \;\;J = \begin{pmatrix} 0 & 1 \\ -1 & 0 \end{pmatrix}.
\]
Since $B = JAJ^{-1}$ it suffices to check that there exist lifts 
$\widetilde{A}, \widetilde{B}, \widetilde{J} \in 
\widetilde{SL}(2,{\mathbb R})$ of $A$, $B$ and $J$ respectively so that
\begin{itemize}
\item $\widetilde{B} = \widetilde{J}\widetilde{A}\widetilde{J}^{-1}$
\item $\widetilde{A}\widetilde{B}\widetilde{A} =
\widetilde{B}\widetilde{A}\widetilde{B}$. 
\item $(\widetilde{A}\widetilde{B})^{6}$ generates 
${\mathbb Z} \cong \ker[\widetilde{SL}(2,{\mathbb Z}) \to SL(2,{\mathbb Z})]$.
\end{itemize}

Let $k \in {\mathbb Z}$. In terms of the description (\ref{usl2})
choose 
\[
\widetilde{A}_{k} = (\begin{pmatrix} 1 & 1 \\ 0 & 1 \end{pmatrix}, 4k)
\quad 
\widetilde{J} = (\begin{pmatrix} 0 & 1 \\ -1 & 0 \end{pmatrix}, 1)
\quad 
\widetilde{B}_{k} = \widetilde{J}\widetilde{A}\widetilde{J}^{-1} =
(\begin{pmatrix} 1 & 0 \\ -1 & 1 \end{pmatrix}, 4k+1).
\]
Using the group law (\ref{grouplaw}) one calculates
\[
\widetilde{A}_{k}\widetilde{B}_{k} = (\begin{pmatrix} 0 & 1 \\ -1 & 1
\end{pmatrix}, 8k + 1 + \tau\left({\mathbb R}p,{\mathbb
R}p,{\mathbb R}q
\right)) = (\begin{pmatrix} 0 & 1 \\ -1 & 1
\end{pmatrix}, 8k + 1).
\]
Similarly we have 
\[
\widetilde{A}_{k}\widetilde{B}_{k}\widetilde{A}_{k} = 
(\begin{pmatrix} 0 & 1 \\ -1 & 0
\end{pmatrix}, 12k + 1 + \tau\left({\mathbb R}p,{\mathbb R}q,{\mathbb R}q
\right)) = (\begin{pmatrix} 0 & 1 \\ -1 & 0
\end{pmatrix}, 12k + 1)
\]
and
\[
\widetilde{B}_{k}\widetilde{A}_{k}\widetilde{B}_{k} = 
(\begin{pmatrix} 0 & 1 \\ -1 & 0
\end{pmatrix}, 12k + 2 + \tau\left({\mathbb R}p,{\mathbb
R}(p-q),{\mathbb 
R}q\right)) = (\begin{pmatrix} 0 & 1 \\ -1 & 0
\end{pmatrix}, 12k + 1).
\]
Thus $\widetilde{A}_{k}\widetilde{B}_{k}\widetilde{A}_{k} =
\widetilde{B}_{k}\widetilde{A}_{k}\widetilde{B}_{k}$ as required. 
Furthermore we
have 
\[
(\widetilde{A}_{k}\widetilde{B}_{k}\widetilde{A}_{k})(\widetilde{B}_{k}
\widetilde{A}_{k}\widetilde{B}_{k}) = (\begin{pmatrix} -1 & 0 \\ 0 & -1
\end{pmatrix}, 24k + 2 + \tau\left({\mathbb R}p,{\mathbb R}q,{\mathbb R}p
\right)) = (\begin{pmatrix} -1 & 0 \\ 0 & -1
\end{pmatrix}, 24k + 2)
\]
and hence
\[
(\widetilde{A}_{k}\widetilde{B}_{k})^{6} = (\begin{pmatrix} 1 & 0 \\ 0 & 1
\end{pmatrix}, 48k + 4).
\]
In addition from the description (\ref{usl2}) we see that 
\[
\ker[\widetilde{SL}(2,{\mathbb Z}) \to SL(2,{\mathbb Z})] = \left\{ (
\begin{pmatrix} 1 & 0 \\ 0 & 1 \end{pmatrix}, 4k) \left| k \in
{\mathbb Z} \right. \right\}
\]
and so $(\widetilde{A}_{k}\widetilde{B}_{k})^{6}$ will generate
$\ker[\widetilde{SL}(2,{\mathbb Z}) \to SL(2,{\mathbb Z})]$ if $k = 0$.
The Lemma is proven. \hfill $\Box$

\begin{rem} \label{rem-angle} 
From the description (\ref{usl2}) it is clear that the
element $\widetilde{A}_{k}$
is the most general lift of $A$. Since any two lifts of $J$ differ by
a central element in $\widetilde{SL}(2,{\mathbb Z})$ it follows that
the pair $(\widetilde{A}_{k}, \widetilde{B}_{k})$ used in the proof of
Lemma~\ref{lem-sl2} represents the most general way to lift $A$ and
$B$ to elements in $\widetilde{SL}(2,{\mathbb Z})$ so that that they
still remain conjugate. The choice of $k = 0$ guaranteeing the
identification of $\op{Map}_{1,1}$ with $\widetilde{SL}(2,{\mathbb
Z})$ has the following simple geometric interpretation.

The natural action of $SL(2,{\mathbb R})$ on the Lagrangian
Grassmanian $\Lambda \cong S^{1}$ lifts to a non-trivial action of the
universal covering group $\widetilde{SL}(2,{\mathbb R})$ on the
universal cover $\widetilde{\Lambda} \cong {\mathbb R}$ of
$\Lambda$. The group ${\mathbb Z} = \ker[\widetilde{SL}(2,{\mathbb R})
\to SL(2,{\mathbb R})]$ acts discretely on $\widetilde{\Lambda}$ with
a fundamental domain isomorphic to an interval of length $2\pi$ and
the action of $k \in {\mathbb Z}$ on $\widetilde{L}$ is by translation
by $2k\pi$. In particular for a
point $\tilde{\ell} \in \widetilde{\Lambda} = {\mathbb R}$ and an
element $\tilde{g} \in \widetilde{SL}(2,{\mathbb R})$ we can
define the {\em displacement angle} of $\tilde{g}$ at $\tilde{\ell}$
as the real number $(\tilde{g}\tilde{\ell} - \tilde{\ell}) \mod
2\pi{\mathbb Z}$. For example if we identify the universal cover
$\widetilde{U}(1) \cong {\mathbb R}$ with a subgroup of
$\widetilde{SL}(2,{\mathbb R})$, then the elements in
$\widetilde{U}(1) \cong {\mathbb R}$ have displacement angles that are
independent of the choice of the point $\tilde{\ell}$ and thus act as
translations on $\widetilde{\Lambda} \cong {\mathbb R}$.

If now $1 \neq g \in SL(2,{\mathbb R})$ is a nilpotent element, then $g$ has
a fixed point on $\Lambda$ and hence among all lifts of $g$ in
$\widetilde{SL}(2,{\mathbb R})$ there is a unique one having a
periodic sequence of fixed points in $\widetilde{\Lambda}$ with period
$\pi$.

The lifts $\widetilde{A}_{0}$ and $\widetilde{B}_{0}$
are characterized uniquely
(see the explanation below) as
the lifts of $A$ and $B$ having periodic sequences of fixed points in
$\widetilde{\Lambda}$. In fact, granted this characterization, one can
easily prove Lemma~\ref{lem-sl2}. Indeed, the matrix $AB \in
SL(2,{\mathbb R})$ is conjugate to a rotation by $-\pi/3$ and so the
displacement angle of any element in $\widetilde{SL}(2,{\mathbb R})$
lifting $AB$ at any point in $\widetilde{\Lambda}$ will be exactly
$-\pi/3$. On the other hand since $\widetilde{A}_{0}$ and
$\widetilde{B}_{0}$ each have a sequence of fixed points which is
periodic with period $\pi$ we see that the displacement angles of
$\widetilde{A}_{0}$ and $\widetilde{B}_{0}$ can not be smaller than
$-\pi$. Thus
both $\widetilde{A}_{0}\widetilde{B}_{0}\widetilde{A}_{0}$ and 
$\widetilde{B}_{0}\widetilde{A}_{0}\widetilde{B}_{0}$ have
displacement angles which are strictly bigger than $-2\pi$ and since
they lift the same element $ABA = BAB \in SL(2,{\mathbb R})$ we must
have $\widetilde{A}_{0}\widetilde{B}_{0}\widetilde{A}_{0} =
\widetilde{B}_{0}\widetilde{A}_{0}\widetilde{B}_{0}$.

To justify the characterization of $\widetilde{A}_{0}$ and
$\widetilde{B}_{0}$ as lifts of $A$ and $B$ having fixed points one
proceeds as follows. Note that similarly to (\ref{usl2}) the space
$\widetilde{\Lambda}$ can be identified
\cite[Section~1.9]{lion-vergne} as a set with
\[
\widetilde{\Lambda} = \{ (\ell,k) \in \Lambda\times {\mathbb Z} | k 
\equiv 1+\dim_{\mathbb R}(\ell\cap \ell_{0}) \mod 2 \}.
\]
Under this identification the action of $(g,n) \in
\widetilde{SL}(2,{\mathbb R})$ on $\widetilde{\Lambda}$ is given by
$(g,n)\cdot (\ell,k) = (g\ell, n + k + \tau(\ell_{0}, g\ell_{0},
g\ell))$ and we immediately see that $\widetilde{A}$ fixes each of the
points $\{({\mathbb R}p, 2k)\}_{k \in {\mathbb Z}}$ and that
$\widetilde{B}$ fixes each of the points $\{ ({\mathbb R}q, 2k+1)
\}_{k \in {\mathbb Z}}$.
\end{rem}

\begin{rem} \label{rem-mumford}  
Lemma~\ref{lem-sl2} (a) easily implies Mumford's 
theorem asserting that $\op{Pic}({\mathcal M}_{1}^{1})$ is a cyclic group
of order 12. Indeed due to Lemma~\ref{lem-sl2} (a) we only need to
check that
the group of central extensions $H^{2}(SL(2,{\mathbb Z}),{\mathbb Z})$
is isomorphic to ${\mathbb Z}/12$ and is generated by the
class of (\ref{meta}). Recall that $SL(2,{\mathbb Z}) = ({\mathbb
Z}/6)\ast_{({\mathbb Z}/2)} ({\mathbb
Z}/4)$ where the cyclic subgroups ${\mathbb
Z}/6$ and ${\mathbb Z}/4$ are generated by the
matrices $AB$ and $ABA$ respectively. The commutator subgroup of
$SL(2,{\mathbb Z})$ is a free subgroup $\Phi$ on two generators,
namely $\Phi$ is the subgroup generated by $[B,A]]$ and $[B,A^{-1}]$.
Hence $SL(2,{\mathbb Z})$ fits in a short exact sequence
$0 \to \Phi \to SL(2,{\mathbb Z}) \to {\mathbb Z}/12 \to 0$
and we have the Hochschild-Serre spectral sequence 
\begin{equation} \label{hs}
E^{pq}_{2} := H^{p}({\mathbb Z}/12,H^{q}(\Phi,{\mathbb
Z}))  \Rightarrow H^{p+q}(SL(2,{\mathbb Z}),{\mathbb Z})
\end{equation}
abutting to the cohomology of $SL(2,{\mathbb Z})$. This is a first
quadrant spectral sequence and since $\Phi$ is of cohomological 
dimension one only the first two rows of (\ref{hs}) 
are non-trivial. In particular (\ref{hs})
degenerates in the $E_{3}$-term and so $E^{02}_{\infty} = E^{02}_{2} =
0$, $E^{11}_{\infty} = \ker d^{11}_{2}$, $E^{20}_{\infty} =
E^{20}_{2}$. Also since ${\mathbb Z}$ is the trivial $SL(2,{\mathbb
Z})$ module we have
\begin{align*}
E^{20}_{2} & = H^{2}({\mathbb Z}/12, H^{0}(\Phi,{\mathbb Z})) = 
H^{2}({\mathbb Z}/12, ({\mathbb Z})^{\Phi})
= H^{2}({\mathbb Z}/12,{\mathbb Z}) \\
E^{11}_{2} & =  H^{1}({\mathbb Z}/12, H^{1}(\Phi,{\mathbb
Z})) =  H^{1}({\mathbb Z}/12, \op{Hom}_{\mathbb Z}(\Phi,{\mathbb
Z})) \\
& = H^{1}({\mathbb Z}/12, {\mathbb Z}^{2}) =
\op{Hom}_{\mathbb Z}({\mathbb Z}/12, {\mathbb Z}^{2}) = 0.
\end{align*}
In  particular we have an isomorphism 
$H^{2}(SL(2,{\mathbb Z}),{\mathbb Z}) = H^{2}(SL(2,{\mathbb
Z})/\Phi,{\mathbb Z}) = H^{2}({\mathbb Z}/12,{\mathbb
Z})$. On the other hand the pullback of the canonical central extension 
\[
0 \to {\mathbb Z} \to {\mathbb Q} \to {\mathbb Q}/{\mathbb Z} \to 0
\]
via any homomorphism ${\mathbb Z}/12 \to {\mathbb
Q}/{\mathbb Z}$ gives an identification $H^{2}({\mathbb Z}/12,{\mathbb
Z})  = \op{Hom}_{{\mathbb Z}}({\mathbb Z}/12,{\mathbb Q}/{\mathbb Z}) 
\cong {\mathbb Z}/12$. 

Finally since both $AB$ and $ABA$ are conjugate to rotations of finite
order we can find maximal compact subgroups $U', U'' \subset
SL(2,{\mathbb R})$ so that $AB \in U'$ and $ABA \in U''$. Furthermore
since the inclusions $U' \subset SL(2,{\mathbb R})$ and $U'' \subset
SL(2,{\mathbb R})$ are homotopy equivalences and since $U'$ and $U''$
are both isomorphic to ${\mathbb R}/{\mathbb Z}$ we see that the pull
backs of the extension (\ref{meta}) by the inclusions ${\mathbb Z}/6
\subset SL(2,{\mathbb Z})$ and ${\mathbb Z}/4 \subset SL(2,{\mathbb
Z})$ are just the standard extensions 
\[
0 \to {\mathbb Z}
\stackrel{6}{\to} {\mathbb Z} \to {\mathbb Z}/6 \to 0 \quad \text{and}
\quad
0 \to {\mathbb Z}
\stackrel{4}{\to} {\mathbb Z} \to {\mathbb Z}/4 \to 0
\]
respectively. In particular the extension (\ref{meta}) is the pullback
of 
\[
0 \to {\mathbb Z}
\stackrel{12}{\to} {\mathbb Z} \to {\mathbb Z}/12 \to 0
\] 
via the
homomorphism $SL(2,{\mathbb Z}) \cong ({\mathbb Z}/6)\ast_{({\mathbb
Z}/2)} ({\mathbb Z}/4) \twoheadrightarrow {\mathbb Z}/12$ and hence
corresponds to the generator of $\op{Hom}_{{\mathbb Z}}({\mathbb
Z}/12,{\mathbb Q}/{\mathbb Z})$. 
\end{rem}

\subsection{Sections in elliptic symplectic fibrations}
\label{subsec-sections}

Now we apply the above considerations to the study of symplectic 
Lefschetz fibrations of genus one.

Consider a symplectic Lefschetz fibration $f : X \to S^{2}$ 
of genus one which corresponds to a relation 
\begin{equation} \label{relation}
t_{1}^{n_{1}}t_{2}^{n_{2}}\ldots
t_{m}^{n_{m}} = 1, \text{ with } n_{i} > 0 \text{ for all } i = 1, \ldots, m.
\end{equation} 
in $\op{Map}_{1}^{1} = SL(2,{\mathbb Z})$ as in
Lemma~\ref{lem-tlf}. To avoid pathologies we will assume that $f$
has a continuous section and that at least one of  the Dehn twists
$t_{i}$ is a twist along a non-separating cycle. In that case it is
known \cite[Theorem~9]{moishezon}
\cite[Theorem~2.4]{kas} that the relation determines $f : X \to
S^{2}$ up to a diffeomorphism. Let $\sigma$ be a smooth section of $f$.
We would like to find the relationship between $\sigma^{2}$ and the
numbers $n_{i}$ (see Corollary~\ref{szpiro}).

The Dehn twists $t_i$ 
correspond to a sequence of smooth circles $s_1, \ldots, s_{m}$ in the generic
fiber $X_o \cong C_{1}$ of $f$. Since the twisting diffeomorphisms
$T_{s_{i}}$ are non trivial only in small neighborhoods of $s_{i}$ in
$C_{1}$ we may assume without loss of generality that there is point
$p\in X_o$  and a small disc $p \in \Delta \subset X_{o}$ so that that all
$T_{s_{i}}$ act identically on $\Delta$. Recall (Lemma~\ref{lem-tlf}) that  
the fibration $f$ was reconstructed from the relation (\ref{relation})
in two steps. We produce first a SLF $u : U \to D$ over a disc $D$ with
a boundary monodromy transformation $\prod_{i}T_{s_{i}}^{n_{i}} : 
C_{1} \to C_{1}$. Next, since by assumption
$\prod_{i}T_{s_{i}}^{n_{i}}$ is homotopic to the identity we can 
complete the resulting fibration into an elliptic fibration over $S^2$
by gluing the trivial genus one family along the boundary of $D$. 

Now we can construct one section where the computation of the
square can be made explicitly. Indeed the fibration $u : U \to D$ has
a section $D\times \{ p \}$ by construction. Also the fiber bundle
$u_{|u^{-1}(\partial D)} :  u^{-1}(\partial D) \to \partial D$ is
homotopy equivalent to $\partial D\times X_{o}$. Since
$\op{Map}_{1}^{1} = \op{Map}_{1}$  we may choose this homotopy
equivalence so that it stabilizes the point $p \in X_{o}$ and thus
glue the section $D\times \{ p \}$ with the constant section trough
$p$ in the trivial family we are gluing to $D$ in order to complete
$u$ to $f$. Let $\sigma$ denote the resulting section.

As we have just seen the right twists $t_{i}$ can be viewed as
elements of the group  $\op{Map}_{1,1}$ which is the group of 
classes of maps stabilizing $\Delta \subset X_{o} = C_{1}$. Then
from the handle body decomposition described in Lemma~\ref{lem-tlf} we
conclude that there is a non-negative integer $n$ so that
\[
\prod t_i^{n_i} = t_{c}^{n}
\]  
where $c$ is the
boundary circle of $C_{1}\setminus \Delta$. In particular the degree
of the normal bundle to $\sigma$ in in $X$ is exactly $-n$.

As a corollary from this observation one can obtain a symplectic 
version of a well known result of 
L. Szpiro \cite{SZ} who proved that for any jacobian 
 elliptic fibration (elliptic fibration with a holomorphic section) 
over ${\mathbb P}^1$ with only multiplicative fibers we have an inequality 
between the number $D$ of singular fibers and the number $N$ of irreducible 
components of the singular fibers. More precisely he showed that $D
\leq 6N$. 

Since in the monodromy group the multiplicative fiber with 
$m$ components corresponds to the element $t^m$ for some right
Dehn twist $t$ we see that the following Corollary is a straightforward
generalization of Szpiro's result to the symplectic category.

\begin{corr} Assume that in $\op{Map}_{1,1} =
\widetilde{SL}(2,{\mathbb Z})$ we have a relation 
$\prod_{i = 1}^{m} t_i^{n_i} = t_{c}^{n}$ with $n_i > 0$. Then
 $\sum n_i = 12n$ and  $m > 2n$. \label{szpiro}
\end {corr}
{\bf Proof.} From the standard presentation (\ref{presentation}) 
of $\op{Map}_{1,1}$ we see
that there exists a unique homomorphism
$\chi : \op{Map}_{1,1} \to {\mathbb
 Z}$ characterized by the property that $\chi(t_{a}) = 1$ and
$\chi(t_{b}) = 1$. Alternatively by Remark~\ref{rem-mumford} the
pushout of the extension (\ref{metaplectic}) via the multiplication
map $\op{mult}_{12} : {\mathbb Z} \to {\mathbb Z}$ is a split
extension and so we can choose for $\chi$ the composition of
$\op{mult}_{12}$ with a splitting of this extension. 
Since $t_{c} = (t_{a}t_{b})^{6}$ we have that
$\chi(t_{c}) = 12$ and that 
\[
\sum n_{i} = \chi(\prod_{i = 1}^{m} t_i^{n_i}) = \chi(t_{c}^{n}) = 12n.
\]
Furthermore as explained in Remark~\ref{rem-angle} 
the transformation $t_{c}^{n}$ translates every point in 
${\mathbb R} = \widetilde{\Lambda}$ by $-2n\pi $. On the other hand 
$t_i^{n_i}$ projects to a nilpotent transformation in $SL(2,{\mathbb
Z})$ and so has a displacement angle $-\pi$ at any point in 
$\widetilde{\Lambda}$. Hence $m \geq 2n$. \hfill $\Box$

\bigskip

We are now in a position to prove Corollary~\ref{identity}

\bigskip

\noindent
{\bf Proof of Corollary~\ref{identity}.} If the subsemigroup in
$\widetilde{SL}(2,{\mathbb Z})$ generated by the conjugates of $t_{a}$
contains the identity element, then we can find a relation of the form
(\ref{relation}) in $\op{Map}_{1,1}$. Hence there exist a SLF $f : X
\to S^{2}$ where the section corresponding to the puncture $p$ has a
trivial tubular neighborhood and normal bundle. By the above
calculation this yields $0 = \sigma^{2} = - \sum n_{i}$, i.e. $n_{i} =
0$ for all $i$. This proves the corollary. \hfill $\Box$

\begin{rem}  It is very tempting  to try  to extend this purely
group theoretic proof of Szpiro's result to the elliptic curves over number 
fields. It is expected that the analogue of Szpiro's inequality
with any constant instead of $6$ will lead to a  solution of the  
ABC-conjecture.
\end{rem}

\subsection{Symplectic fibrations with fibers of higher genus}

Some of the features of $\op{Map}_{1,1}$ that allowed us to restrict
the numerical properties of the sections in a genus one SLF carry over
to the case of SLF of higher genus. Here we outline an approach to the
study of the numerical behavior of the sections in a higher genus SLF.

The main ingredient in the discussion in the previous section was the
existence of a natural action of $\op{Map}_{1,1}$ on the universal
cover $\widetilde{\Lambda} = {\mathbb R}$ of the Lagrangian
Grassmanian $\Lambda$ and the notion of a displacement angle.
 
There is a similar action of $\op{Map}_{g,1}$ which we proceed to
describe. For more details the reader may wish to consult the
excellent exposition of S. Morita \cite{morita}.

As explained in Remark~\ref{remark-extension} the group
$\op{Map}_{g,1}$ is a central extension 
\begin{equation} \label{euler}
0 \to {\mathbb Z} \to \op{Map}_{g,1} \to \op{Map}_{g}^{1} \to 1
\end{equation}
of $\op{Map}_{g}^{1}$ by an infinite cyclic group generated by the Dehn
twist around the puncture.

This central extension has many remarkable properties but we want
to emphasize the ones which are reflected in the structure of the
semi group $\op{Map}_{g,1}(-)\subset \op{Map}_{g,1}$ generated by all
right handed Dehn twists.
 
Let $e \in H^2(\op{Map}_{g}^{1}),{\mathbb Z})$ denote the class of the
central extension (\ref{euler}). The class $e$ is dubbed the {\em
Euler class} by Morita and admits the following simple 
description \cite{morita}.
As in the beginning of Section~\ref{subsec-positive} we will view the
elements in $\op{Map}_{g}^{1})$ as isotopy classes of orientation
preserving diffeomorphisms of $C_{g}$ that preserve the point $x_{1}
\in C_{g}$.

Fix an isomorphism of the universal covering of $C_{g}$ with the unit
disk $D \subset {\mathbb C}$, e.g. by fixing a point $j$ in the
Teichm\"{u}ller space ${\mathcal T}_{g}$. Let $\nu : D \to C_{g}$ be
the corresponding covering map. Any orientation preserving 
diffeomorphism
$\Phi : C_g \to C_g$ with $\Phi(x_{1}) = x_{1}$ defines a quasi-conformal 
map $h^{\Phi} : D \to D$ which preserves the preimage $\nu^{-1}(x_{1})$ 
of $x_{1}$ in $D$.
This map extends to an orientation preserving homeomorphism 
$h^{\Phi}_{\partial} : S^1 \to S^1$ of the boundary $\partial D \cong S^{1}$
of $D$. Modulo isometries of $D$ the homeomorphism
$h^{\Phi}_{\partial}$ is uniquely defined by the isotopy class of 
$\Phi$ relative to $p$. If in addition we require that the isomorphism
$j : \widetilde{C}_{g} \widetilde{\to} D$ sends a marked point
$\tilde{x}_{1} \in \nu^{-1}(x_{1})$ to $0 \in D$ and the differential
$dj$ induces a fixed isomorphism of the tangent spaces
$T_{\tilde{x}_{1}}\widetilde{C}_{g}$ and $T_{0}D$, then the
element $h^{\Phi}_{\partial} \in \op{Homeo}^{+}(S^{1})$ 
is uniquely determined by the mapping class $[\Phi] \in
\op{Map}_{g}^{1}$ of $\Phi$.

Therefore we obtain a natural action $\rho : \op{Map}_{g}^{1} \to
\op{Homeo}^{+}(S^{1})$ given by $\rho([\Phi]) = h^{\Phi}_{\partial}$. 
The group $\op{Homeo}^{+}(S^1)$ has a natural central extension by an
infinite cyclic group which is the group
$\op{Homeo}^{\op{per}}({\mathbb R})$
of $2\pi$-periodic
orientation preserving 
homeomorphisms  of the line
${\mathbb R}$ (here as before one should view the line ${\mathbb R}$
as the universal cover of the circle $S^{1}$). One can check
\cite{morita-characteristic}
that the central extension (\ref{euler}) is just the pullback of the 
extension $0 \to {\mathbb Z} \to \op{Homeo}^{\op{per}}({\mathbb R})
\to \op{Homeo}^{+}(S^1) \to 1$ by the homomorphism $\rho$.
In particular we have a natural homomorphism $\rho^{\op{per}} : 
\op{Map}_{g,1} \to \op{Homeo}^{\op{per}}({\mathbb R})$ and so in the
same way as in Remark~\ref{rem-angle} we can
define the {\em displacement angle} associated with any pair
$(\phi,\ell)$, $\phi \in $, $\ell \in {\mathbb R}$ as the real number 
$(\rho^{per}(\phi)(\ell) - \ell) \mod 2\pi{\mathbb Z}$.

We propose the following conjectural characterization of the semigroup
$\op{Map}_{g,1}(-)$ which generalizes Remark~\ref{rem-angle}: 

\begin{con} \label{con-genus-g-semigroup}
The subsemigroup $\op{Map}_{g,1}(-) \subset \op{Map}_{g,1}$ consists
of mapping classes whose displacement angle is non-negative at
every point in ${\mathbb R}$.
\end{con}

\noindent
Before we give some evidence for the validity of this conjecture we
need to recall Morita's analysis of the Euler class $e$.

\medskip

Similarly to the genus one case we can view the Euler class $e$ as 
an element of Picard group of the
moduli space ${\mathcal M}_{g}^{1}$ of smooth curves of genus $g$ 
with one marked point. As such $e$ can be written explicitly
\cite{harris}, \cite{morita-familiesii} 
as a linear combination with ${\mathbb Q}$-coefficients of two
cocycles of geometric origin: $e_{1}$ and $c$.  
The cocycle $e_{1}$ is
is just the first Mumford-Miller-Morita class and is proportional
to the generator of the Picard group of the moduli space ${\mathcal
M}_{g}$ of smooth curves of genus $g$. Also as in the genus one case
the cocycle $e_{1}$ is induced from a cocycle for
$Sp(H_{1}(C_{g},{\mathbb Z}))$
via the natural homomorphism $\op{Map}_{g}^{1} \to \op{Map}_{g} \to
Sp(H_{1}(C_{g},{\mathbb Z}))$. 
The cocycle $c$ is somewhat more mysterious. Geometrically it can be
defined as follows. Let $\Theta$ denote the pullback of the first
Chern class of the relative theta line bundle on the universal degree
$g - 1$ Jacobian ${\mathcal J}^{g-1} \to {\mathcal M}_{g}$ via the
canonical Abel-Jacobi map from the universal curve ${\mathcal M}_{g}^{1}
\to {\mathcal M}_{g}$ to ${\mathcal J}^{g-1} \to {\mathcal M}_{g}$
defined by the point $x_{1} : {\mathcal M}_{g} \to {\mathcal
M}_{g}^{1}$. Then it can be shown that $c = 8\Theta - e_{1}/3$
\cite[(1.7)]{morita-familiesii} and \cite[Proposition~2]{hain-reed}.  

Observe that these geometric definitions of the classes $e_{1}$ and
$c$ make sense in the genus one case as well. However it follows from
the proof of Mumford's result \cite[Main~Theorem]{mumford-picard} discussed
in Section~\ref{subsec-mapone} that in the genus one case the two
classes $e_{1}$ and $c$ coincide and so we are not getting anything
new for the group $\op{Map}_{1} = SL(2,{\mathbb Z})$.
In contrast for $g > 2$ the classes $e_{1}$ and $c$ are linearly
independent in $H^{2}(\op{Map}_{g}^{1},{\mathbb Q})$ and are related
to the class $e$ by the following explicit formula (see
e.g. \cite{morita-familiesii}) 
\[
e = \frac{1}{4g(1-g)}(e_{1} + c).
\]
Algebraically the classes $e_{1}$ and $c$ can be interpreted as
follows. Consider the lattice $H := H_{1}(C_{g}, {\mathbb Z})$
together with its natural symplectic form $\theta$ 
given by the intersection pairing. Morita has shown (see e.g. 
\cite{morita-characteristic}) that there exists a two-step nilpotent group
${\mathcal U}$ equipped with a $Sp(H,\theta)$-action and a homomorphism
$\kappa : \op{Map}_{g}^{1} \to Sp(H,\theta)\ltimes ({\mathcal U}\otimes
{\mathbb Q})$ so that the classes $e_{1}$ and $c$ are pullbacks of
natural cohomology classes on $Sp(H,\theta)\ltimes ({\mathcal U}\otimes
{\mathbb Q})$. More precisely ${\mathcal U}$ 
is a $Sp(H,\theta)$-invariant central extension of
$\wedge^{3}H$ by a certain $Sp(H,\theta)$-module  of
finite rank and $e_{1}$ is a pullback of a cohomology class of $Sp(H,\theta)$
and $c$ is a pullback of a cohomology class $e_{1}$ of $Sp(H,\theta)\ltimes
\frac{1}{2}\wedge^{3}H$ via the natural maps 
\begin{align}
\op{Map}_{g}^{1} \stackrel{\kappa}{\to} & Sp(H,\theta)\ltimes 
({\mathcal U}\otimes {\mathbb Q}) \to  Sp(H,\theta) \label{quot1}\\
\op{Map}_{g}^{1} \stackrel{\kappa}{\to} & Sp(H,\theta)\ltimes 
({\mathcal U}\otimes {\mathbb Q}) \to Sp(H,\theta)\ltimes
(\wedge^{3}H\otimes {\mathbb Q}) \label{quot2}
\end{align}
respectively.

To identify the element in $H^{2}(Sp(H,\theta),{\mathbb
Z})$ that pulls back to $e_{1}$ consider the $2g$-dimensional vector
space $H_{{\mathbb R}} := H\otimes
{\mathbb R}$.
The symplectic group $Sp(H_{{\mathbb R}},\theta)$ is homotopy
equivalent to its maximal compact subgroup which in turn is isomorphic
to the unitary group $U(g)$. In particular $\pi_{1}(Sp(H_{{\mathbb
R}},\theta)) \cong \pi_{1}(U(g)) \cong {\mathbb Z}$ and so the
universal cover $\widetilde{Sp}(H_{{\mathbb R}},\theta)$ of 
$Sp(H_{{\mathbb R}},\theta)$ is naturally a  central extension of 
$Sp(H_{{\mathbb R}},\theta)$ by an infinite cyclic group. By pulling
back this extension by the natural inclusion $Sp(H,\theta) \subset 
Sp(H_{{\mathbb R}},\theta)$ we get a central extension
\begin{equation} \label{meta2g}
0 \to {\mathbb Z} \to \widetilde{Sp}(H,\theta) \to
Sp(H,\theta) \to 1.
\end{equation}
From the geometric description of $e_{1}$ given above it is clear that
$e_{1}$ is proportional to the pullback of the extension class
(\ref{meta2g})\footnote{In fact it was pointed out to us by R.~Hain that
the extension class (\ref{meta2g}) determines the universal central
extension of $Sp(H,\theta)$ and so with the correct choice of the
orientation of the fiber must be equal to the pullback of the
determinant of the pushforward of the relative cotangent bundle on the
universal abelian variety or in other words to $e_{1}/12$.}.

To identify the element in $H^{2}(Sp(H,\theta)\ltimes
\frac{1}{2}\wedge^{3}H, {\mathbb Z})$ that pulls back to $c$ notice
that the contraction with $\theta$ gives a well defined homomorphism
of $Sp(H,\theta)$-modules $C : \wedge^{3}H \to H$. Consider the
Heisenberg central extension $0 \to {\mathbb Z} \to {\mathcal H} \to H
\to 0$ corresponding to the class $\theta \in H^{2}(H,{\mathbb
Z})$. Since by definition $\theta$ is $Sp(H,\theta)$-invariant the
pullback of the Heisenberg extension via the map $C$ will also be 
$Sp(H,\theta)$-invariant and will so determine an element 
$H^{2}(Sp(H,\theta)\ltimes
\frac{1}{2}\wedge^{3}H, {\mathbb Z})$. In 
\cite[Theorem~3.1]{morita-characteristic} Morita shows 
that this element pulls back
to $c$ via the homomorphism (\ref{quot2}).

It is not hard to see that the above interpretation of $c$  gives
rise to a representation of $\op{Map}_{g,1}$ into the universal cover
of a different symplectic
group. Indeed, let us fix a symplectic basis of $H$ and  identify
$Sp(H,\theta)$ with the group $Sp(2g, {\mathbb Z})$. 
Denote by $ASp(2g,{\mathbb Z})$ the subgroup of
the group $Sp(2g + 2, {\mathbb Z})$ which stabilizes the orthogonal
complement of the element $b_{g+1}$ of the standard symplectic basis $a_{1},
b_{1}, \ldots, a_{g+1}, b_{g+1}$ of ${\mathbb Z}^{2g + 2}$. 
The natural map $ASp(2g, {\mathbb Z}) \to GL(2g+1,{\mathbb Z})$, $X
\mapsto  X_{|b_{g+1}^{\perp_{\theta}}}$ maps the group $ASp(2g, {\mathbb
Z})$ onto $Sp(2g,{\mathbb Z})\ltimes {\mathbb Z}^{2g}$ and so 
$ASp(2g, {\mathbb Z})$ fits into a central extension
\begin{equation}
0 \to {\mathbb Z} \to ASp(2g, {\mathbb
Z}) \to Sp(2g,{\mathbb Z})\ltimes {\mathbb Z}^{2g} \to 1
\end{equation}
which is clearly isomorphic with the central extension of
$Sp(H,\theta)\ltimes H$ considered above.

The fact that $e$ is proportional to $e_{1} + c$ combined with the
above discussion implies that the group $\op{Map}_{g,1}$ has a 
well defined homomorphism
$\delta : \op{Map}_{g,1} \to \widetilde{ASp}(2g, {\mathbb Z})$ to the universal
cover
$\widetilde{ASp}(2g, {\mathbb Z})$  of $ASp(2g, {\mathbb Z})$. 

On the other hand it is clear from the construction that
$\widetilde{ASp}(2g, {\mathbb Z})$ is also isomorphic to the preimage
of $ASp(2g, {\mathbb Z})$ into the
universal cover of $Sp(2g+2, {\mathbb R})$. The choice of an inclusion 
$U(g+2) \subset Sp(2g+2,{\mathbb R})$ induces an isomorphism of 
the Lagrangian Grassmanian $\Lambda$ of ${\mathbb R}^{2g+2}$  
with the quotient $U(g+1)/O(2g+2)$. Since $O(2g+2) \subset
SU(g+1)$ we have a well defined map $\det : \Lambda \to U(g+1)/SU(g+1) \cong
S^{1}$ and in particular a well defined map of the universal covers
$\widetilde{\det} : \widetilde{\Lambda} \to {\mathbb R}$.
Since the group $\widetilde{Sp}(2g+2,{\mathbb R})$ acts
naturally on $\widetilde{\Lambda}$ we can use the map 
$\widetilde{\det}$ to define a displacement angle for any element $A \in
\widetilde{Sp}(2g+2,{\mathbb R})$  and any $\ell \in {\mathbb R}$. 

Since the homomorphism $\delta : \op{Map}_{g,1} \to \widetilde{ASp}(2g,
{\mathbb Z})$ is essentially given by the class $e$ it is reasonable
to expect that for any $\phi$ and any $\ell \in {\mathbb R}$ the
displacement angles of $(\rho^{\op{per}}(\phi),\ell)$ and
$(\delta(\phi),\ell)$ will have the same sign. If this is the case, 
then the validity of Conjecture~\ref{con-genus-g-semigroup} 
will easily follow since
similarly to the genus one case we can split the all unipotent
elements in $ASp(2g, {\mathbb Z})$ into two classes 
according to their displacement angle. That is  - unipotent elements
whose lifts in $\widetilde{ASp}(2g, {\mathbb Z})$ have a periodic
sequence of points with displacement angle zero and unipotent elements
elements whose lifts in $\widetilde{ASp}(2g, {\mathbb Z})$ do not have
points with displacement angle zero. 

Since the Dehn twists obviously belong to the first class and since
the elements in the first class will generate the subsemigroup in
$\widetilde{ASp}(2g, {\mathbb Z})$ consisting of elements with
non-positive displacement angle this argument will prove
Conjecture~\ref{con-genus-g-semigroup}.

Similar arguments should also lead to a proof of the following
theorem, which  is obviously correct in the case of projective
Lefschetz pencils.

\begin{theo}[I.Smith] \label{theo-torelli} There are no non-trivial 
SLF whose monodromy group is contained in the Torelli group.
\end{theo}

A special case of the above theorem was originally 
proven by B.Ozbagci 
\cite[Corollary~7]{ozbagci-signatures} who showed that 
hyperelliptic SLF (and in particular  all Lefschetz fibrations of 
fiber genus two) cannot have monodromy contained in the Torelli group. 
The general statement of Theorem~\ref{theo-torelli} 
appeared as a conjecture in a preliminary version of the present paper.
The conjecture was settled affirmatively by Ivan Smith who graciously
provided us with the proof appearing in
Appendix~\ref{appendix-torelli} below.

\

\bigskip

\noindent
In this direction we would like to ask a couple of questions:

Let ${\mathcal L} \to {\mathcal M}_{g}$ be the Hodge line bundle and let
$c_{1}({\mathcal L}) \in CH^{1}(\overline{{\mathcal M}}_{g},{\mathbb
Q})$ denote
the natural extension of the first Chern class of ${\mathcal L}$ to the
Deligne-Mumford compactification of ${\mathcal M}_{g}$.

\begin{que} Let $f : X \to C$ be a symplectic Lefschetz fibration of
fiber genus $g$ over an arbitrary Riemann surface $C$. Is it true that
\[
\langle c_1({\mathcal L}),[C]\rangle \geq  0? 
\]
\end{que}

\

\medskip
This question has an affermative answer (see
Corollary~\ref{ivan-inequality}) whenever $C$ is of genus zero.

\begin{que} Let $f : X \to C$ be a symplectic Lefschetz fibration of
fiber genus $g$ over an arbitrary Riemann surface $C$. Is it true that
\[
(8g+4)\langle c_1({\mathcal L}),[C]\rangle  - g\cdot\mu \geq  0? 
\]
\end{que}

\

\noindent
Here as usual $\mu$ denotes the number of the singular fibers in 
the fibration $f: X \to C$.

The last question is a symplectic version of the Moriwaki inequality
and was suggested by R. Hain.

If the answers of the above questions are positive we get additional
restrictions on the words in the mapping class group defining
SLF. Some analogues of  Corollary~\ref{szpiro} for high genus SLF can be
expected as well.
\bigskip

\

\newpage

\appendix

\Appendix{(by Ivan Smith) Torelli fibrations} 
\label{appendix-torelli}

\

\bigskip

\noindent
The purpose of this Appendix is to present a short proof of
Theorem~\ref{theo-torelli}.

Suppose $f : X \to S^{2}$ is a SLF whose monodromy is contained in the
Torelli group. The SLF $f$ induces a sphere $S \cong S^{2}$  in the 
compactified moduli space of genus $g$  curves.  Indeed let us choose
an almost complex structure on $X$ compatible with the symplectic 
form. The restriction of the almost complex structure on each smooth
fiber is integrable and so one obtains a map $u$ of a punctured sphere
into ${\mathcal M}_g$. By
assumption the map $f$ has a local complex model near the
singular points, and it is easy to see that this gives an integrable almost
complex structure in an entire neighborhood of the singular
fibers. Thus we
can extend $u$ smoothly to a map of the closed sphere into the
compactified moduli space.  The isotopy class of this map is independent of
the choices of almost complex structures on $X$ and in the
neighborhoods of the singular fibers. 

Let ${\mathcal L} \to {\mathcal M}_{g}$ be the Hodge line bundle and let
$c_{1}({\mathcal L}) \in CH^{1}(\overline{{\mathcal M}}_{g},{\mathbb
Q})$ denote
the natural extension of the first Chern class of ${\mathcal L}$ to the
Deligne-Mumford compactification of ${\mathcal M}_{g}$. We will need
the following two preliminary lemmas:

\begin{lem}  \label{lem-ivan-formula} 
For any symplectic Lefschetz fibration, 
\[
\op{sign}(X) = 4\langle c_1({\mathcal L}), [S] \rangle - \mu,
\] 
where $\mu$ is the number of the
singular fibers of $f$ and $\op{sign}(X)$ is the signature of $X$.
\end{lem}
{\bf Proof.} See \cite{smith}. \hfill $\Box$

\begin{lem}  For a fibration with only separating 
vanishing cycles, $\op{sign}(X) =-\mu.$
\end{lem}
{\bf Proof.} See \cite{ozbagci-signatures}. \hfill $\Box$

\bigskip

\noindent
As a consequence of the above two lemmas we get that 
for a symplectic Lefschetz fibration with monodromy group
contained in the Torelli group:
\[
\langle c_1({\mathcal L}), [S] \rangle =  0. 
\]
Therefore Theorem~\ref{theo-torelli} will be proven if we know the following:

\begin{lem} Let $X$ be a symplectic Lefschetz fibration with monodromy
group contained in the Torelli group. Then:
\[
\langle c_1({\mathcal L}), [S] \rangle  >  0. 
\]
\end{lem}
\begin{rem} We need base $S^2$  for this
entire argument. There are fibrations over $T^2$ with monodromy
group in the Torelli group and no singular fibres.
\end{rem}
{\bf Proof} Let $n$ denote the number of exceptional $(-1)$-spheres in
the space  $X$ and let $\varepsilon : X \to 
X_{\op{min}}$ denote the contraction of  all the $(-1)$ spheres.
Clearly since each new exceptional sphere contributes
a homology class, we know  that $n$ is bounded above by the second
Betti number $b_2(X)$. In fact without a
loss of generality we may assume that $n \leq b_{2}(X) - 2$. Indeed
the fibers of $f : X \to S^{2}$ are not exceptional, and we can fiber
sum  $X$  with itself to
ensure that $f : X \to S^{2}$ has a section with even square and in 
particular is not exceptional.

Furthermore we may assume that $X_{\op{min}}$ is not symplectomorphic
to an irrational ruled surface. Indeed if $\pi : X_{\op{min}} \to C$
is a sphere bundle over a surface $C$ of genus $\geq 1$, then
$H_{1}(C,{\mathbb Z}) = H_{1}(X_{\op{min}},{\mathbb Z}) =
H_{1}(X,{\mathbb Z})$. On the other hand since the geometric 
monodromy of $f$ is contained in the Torelli group it follows that 
$H_{1}(X_{s},{\mathbb Z}) = H_{1}(X,{\mathbb Z})$ for all smooth
fibers $X_{s} \subset X$ of $f$. In particular $g(C) = g(X_{s}) = g$
and so $\pi_{|X_{s}} : X_{s} \to C$ must be a symplectomorphism. This
shows that the natural map $(\pi\circ \varepsilon, f) : X \to C\times
{S^{2}}$ is a symplectic isomorphism in the complement of all $(-1)$
spheres. This however implies that $f$ is a trivial fibration which is
a case we exclude.

Since $X$ is a symplectic Lefschetz fibration of
fiber genus $g$ one has $c_2 (X) = 4 - 4g + \mu$.
Also observe that $c_2(X) = \op{euler}(X) = 2 - 2b_1(X) + b_2(X)$
and   $b_1(X) = 2g$ since all vanishing
cycles are by assumption null-homologous.  So we get the estimate:
\[ 
b_2(X) = 2 + \mu,
\]
and hence 
\[   
n \leq b_2(X)-2  = \mu. 
\]
Now    $c_1^2 (X) + c_2 (X) = c_1^2(X_{\op{min}}) - n + c_2(X)$ and
using  the above estimate on $n$, we see that
\[
c_1 ^2 (X) + c_2 (X) > c_1^2 (X_{\op{min}}) + 2 - 4g \geq  2 - 4g
- 2,
\]
The last inequality follows from a powerful theorem of A.K.Liu
\cite{liu} asserting that for any {\em minimal} symplectic four-manifold $Y$
which is not irrational ruled one has  $c_1^2(Y) \geq 0$. 

On the other hand the signature formula in
Lemma~\ref{lem-ivan-formula} gives
\begin{align*} 
\frac{1}{12}[c_1 ^2(X) + c_2 (X)] & = \frac{1}{4} 
[\op{sign}(X) + c_{2}(X)] \\ 
& =  \frac{1}{4}[4\langle c_1({\mathcal
L}),[S]\rangle - \mu + 4(1-g) + \mu] \\ & \
= \langle c_1({\mathcal L}),[S]\rangle + 1- g
\end{align*}
and hence   
\[
\langle c_1({\mathcal L}),[S]\rangle > \frac{4g - 5}{6},
\]
i.e. $\langle c_1({\mathcal L}),[S]\rangle > 0$ as long as $g >
1$. This proves the lemma since the case $g = 1$ is clear.  \hfill $\Box$

\bigskip

As a consequence of Theorem~\ref{theo-torelli} we get

\begin{lem} Let $X$ be an arbitrary symplectic Lefshetz fibration. Then
\[
\op{sign}(X) +  \mu > 0.  
\]
\end{lem}
{\bf Proof} As it follows from Theorem~\ref{theo-torelli} 
there always exists a
non homologous to zero vanishing cycle in $X$. By applying the local
signature formula of \cite{ozbagci-signatures} we see that the contribution of 
this singular fiber to $\op{sign}(X)$ is either zero or one.
The latter implies the lemma. \hfill $\Box$

As corollary we get:

\begin{corr} \label{ivan-inequality} 
Let $X$ be an arbitrary symplectic Lefshetz fibration. Then:
\[
\langle c_1({\mathcal L}),[S]\rangle >  0. 
\]
\end{corr}

\begin{rem} The above lemmas puts many restrictions  on the
possible monodromies for hyperelliptic families (see
\cite{ozbagci-signatures}).
\end{rem}

\begin{rem} The proof of Theorem~\ref{theo-torelli} is related to 
a question of Gompf who asked whether $c_2(X)$
is positive for symplectic  four-manifolds which are not irrational ruled.
In particular it is an interesting question if the minimal number
of the singular fibers in a SLF is at least $4(g-1)$ if $X$ is  not
irrational ruled symplectic  four-manifold.
\end{rem}

\bigskip

\noindent
{\sc Ivan Smith, New College, Oxford University, \\
smithi@maths.ox.ac.uk}

{\sc
\noindent
J.Amor\'os, UPC, Barcelona, \\amoros@ma1.upc.es 

\bigskip

\noindent
F.Bogomolov, Courant Institute, NYU  \\ bogomolo@MATH8.CIMS.NYU.EDU

\bigskip

\noindent
L.Katzarkov, UC Irvine. \\ lkatzark@math.uci.edu

\bigskip

\noindent
T. Pantev, UPenn, \\
tpantev@math.upenn.edu 
}

\end{document}